\documentclass[10pt]{amsart}
\usepackage{amssymb,amsfonts,amsmath,amsthm}
\usepackage[all]{xy}

\newcommand{\proj}{\mathbb P}

\newcommand{\aff}{\mathbb A}
\newcommand{\bee}{\mathbb B}

\newcommand{\kay}{{\mathbb K}}
\newcommand{\com}{\mathbb{C}}
\newcommand{\zee}{{\mathbb{Z}}}
\newcommand{\zeeg}{{\zee}_{\geq 0}}

\newcommand{\oh}{{\mathcal{O}}}
\newcommand{\ca}{{\mathcal{A}}}

\newcommand{\ci}{{\mathcal{I}}}

\newcommand{\cl}{{\mathcal{L}}}
\newcommand{\cz}{{\mathcal{Z}}}

\newcommand{\cw}{{\mathcal{W}}}
\newcommand{\ex}{{\mathcal{X}}}

\newcommand{\yi}{{\mathcal{Y}}}

\newcommand{\cm}{{\mathcal{M}}}
\newcommand{\cp}{{\mathcal{P}}}

\newcommand{\ma}{{\mathcal{MA}}}
\newcommand{\mz}{{\mathcal{MZ}}}

\newcommand{\mt}{{\mathcal{MT}}}

\newcommand{\acm}{{\mathfrak{M}}}

\newcommand{\Top}{\operatorname{L^0}}
\newcommand{\Bot}{\operatorname{L^\infty}}
\newcommand{\Split}{\operatorname{Split}}
\newcommand{\Dil}{{\operatorname{Dil}}}
\def\Lnt#1{\operatorname{L}_{#1,\text{not top}}}
\def\Lnb#1{\operatorname{L}_{#1,\text{not bot}}}
\def\Le#1{\operatorname{L}_{#1,\text{ext}}}

\newcommand{\ev}{\text{ev}}
\newcommand{\Ev}{\text{Ev}}
\newcommand{\Bl}{\text{Bl}}

\newcommand{\id}{\text{id}}

\newcommand{\TC}{\operatorname{TC}}
\newcommand{\Spec}{\operatorname{Spec}}
\newcommand{\Proj}{\operatorname{Proj}}
\newcommand{\Sing}{\operatorname{Sing}}

\newcommand{\rel}{\text{rel}}
\newcommand{\Obj}{\operatorname{Obj}}
\newcommand{\Mor}{\operatorname{Mor}}

\newcommand{\Sch}{\text{(Sch)}}

\newcommand{\Shom}{\mathcal{H}\text{om}}

\def\sig#1{\sigma_{#1}}
\def\vir#1{[#1]^{\text{vir}}}
\def\Gt#1{\widetilde{G[#1]}}

\newcommand{\Aut}{\operatorname{Aut}}

\def\nicechartonZ#1#2
{\xymatrix{
\ex \ar[d] \ar[r]^{#1} & Z[n]\ar[d]\\
S \ar[r]_{#2} & \aff^n}}
\def\nicechartonZmp#1#2
{\xymatrix{
\ex \ar[d] \ar[r]^{#1} & Z[n]\ar[d]\\
S \ar@<1ex>[u]^{\gamma}\ar[r]_{#2} & \aff^n}}
\def\nicechartonAmp#1#2
{\xymatrix{
\ex \ar[d] \ar[r]^{#1} & A[n]\ar[d]\\
S \ar@<1ex>[u]^{\gamma}\ar[r]_{#2} & \aff^n}}

\newcommand{\srarr}{\rightarrow}

\newcommand{\cs}{\com^*}

\theoremstyle{plain}
\newtheorem{theorem}{Theorem}
\newtheorem{lemma}[theorem]{Lemma}
\newtheorem{proposition}[theorem]{Proposition}
\newtheorem{corollary}[theorem]{Corollary}
\theoremstyle{definition}
\newtheorem{definition}[theorem]{Definition}

\newtheorem{claim}[theorem]{Claim}

\numberwithin{theorem}{subsection}

\begin{document}

\title{Line-Bundles on Stacks of Relative Maps}
\author{Eric Katz}
\date{\today}
\address{Department of Mathematics, Duke University}\email{eekatz@math.duke.edu}
\begin{abstract}
We study line-bundles on the moduli stacks of relative stable and
rubber maps that are used to define relative Gromov-Witten
invariants in the sense of J. Li \cite{Li1,Li2}.  Relations
between these line-bundles yield degeneration formulae which are
used in \cite{Ka2}.  In addition we prove the {\em Trivial
Cylinders Theorem}, a technical result also needed in \cite{Ka2}.
\end{abstract}
\maketitle

    \section{Introduction}

Relative Gromov-Witten invariants \cite{Li1,Li2} are intersections
numbers on a stack parameterizing stable maps to a projective
manifold $Z$ relative to a smooth divisor $D$.  One looks at
stable maps to $Z$ where all points of intersections of the map
with $D$ are marked and multiplicities at these points are
specified.  To obtain a proper moduli stack of such maps, one must
allow the target to degenerate to $\ _k Z=Z\sqcup_D P_1 \sqcup_D
\dots \sqcup_D P_k$, that is, $Z$ union a number of copies of
$P=\proj_D (N_{D/Z}\oplus 1_D)$ the projective completion of the
normal bundle to $D$ in $Z$. Maps with a non-smooth target are
said to be {\em split maps}.  Li constructed a moduli stack of
relative maps called $\cm(\cz,\Gamma)$ for $\Gamma$, a certain
kind of graph, and constructed its virtual fundamental cycle. This
stack has an evaluation map
$$\Ev_{\mz}:\cm(\cz,\Gamma)\srarr Z^m \times D^r$$
where $m$ and $r$ are the number of interior and boundary marked
points, respectively. {\em Relative Gromov-Witten invariants} are
given by evaluating pullbacks of cohomology classes by $\Ev$
against the virtual cycle.

It is natural to break the target $\ _k Z$ as the union of $\ _l
Z=Z\sqcup_D P_1 \sqcup_D \dots \sqcup_D P_l$ and $\ _{k-l-1}
P=P_{l+1}\sqcup_D \dots \sqcup_D P_k$. In fact, such splitting is
necessary to parameterize fixed loci in $\com^*$-localization in
the sense of \cite{K} and \cite{GP} in the relative framework
\cite{GV}. If we set $X=D$, and $L=N_{D/Z}$, the normal bundle to
$D$ in $Z$, one is led to study stable maps into the
projectivization of a line bundle $P=\proj_X(L\oplus 1_X)$
relative to the zero and infinity sections, $D_0$ and $D_\infty$
where two stable maps are declared equivalent if they can be
related by a $\com^*$-factor dilating the fibers of $P\srarr X$.
One can construct a moduli stack of such maps, $\cm(\ca,\Gamma)$
and its virtual cycle.  This moduli stack has certain natural line
bundles, called the target cotangent line bundles, $\Top$ and
$\Bot$ and has an evaluation map
$$\Ev_{\ma}:\cm(\ca,\Gamma)\srarr X^m \times X^{r_0} \times X^{r_\infty}$$
The {\em rubber invariants} are obtained by evaluating pullbacks
of cohomology by $\Ev$ map and powers of $c_1(\Bot)$ against the
virtual cycle.

One can take fiber products of stacks $\cm(\cz,\Gamma_Z)$ and
$\cm(\ca,\Gamma_A)$ to parameterize split maps in a bigger moduli
stack $\cm(\cz,\Gamma_Z*\Gamma_A)$.  Such split maps form divisors
that whose fundamental classes are (up to multiplicity) the first
Chern classes of particular line-bundles on
$\cm(\cz,\Gamma_Z*\Gamma_A)$.  By finding relations between such
line-bundles, one is able to prove {\em degeneration formulae}
that allow us to rewrite Gromov-Witten invariants of one moduli
stack in terms of Gromov-Witten invariants are smaller moduli
stacks.  The approach of this paper is to do an in-depth study of
natural line-bundles on stacks.

In section 2, we introduce some technical definitions.

In section 3, we review the construction of $\cm(\cz,\Gamma_Z)$
and $\cm(\ca,\Gamma_A)$.  This material is essentially a
rephrasing of parts of \cite{Li1,Li2}.

In section 4, we introduce line-bundle systems.  These are
line-bundles defined on a sequence of spaces that obey certain
transition properties under inclusions and a group action.  In
more esoteric language, line-bundle systems are line-bundles on
the stacks of degenerations $\cz^\rel$, $\ca^\rel$ and their
universal targets.

In section 5, we explain how line-bundle systems naturally induce
line-bundles on $\cm(\cz,\Gamma_Z)$ and $\cm(\ca,\Gamma_A)$.  We
describe line-bundles on $\cm(\cz,\Gamma)$: $\Dil$ and $\Le{i}$;
and line-bundles on $\cm(\ca,\Gamma)$: $\Split$, $\Top$, $\Bot$,
$\Lnt{i}$, $\Lnb{i}$.  These line-bundles have geometric meaning:
$\Le{i}$ is a line-bundle which has a section whose zero-stack
consists of maps $f:C\srarr\ _k Z$ so that the $i$th marked point
is not mapped to $Z\subset\ _k Z$ (counted with multiplicity);
$\Split$ is a line-bundle whose zero-stack is all split maps;
$\Lnt{i}$, where $i$ is the label of interior marked point, is a
line-bundle whose zero-stack consists of all split maps
$f:C\srarr\ _k P$ where $i$th marked point is not mapped to $P_k$;
$\Lnb{i}$ is its upside-down analog. These line-bundles satisfy
certain relations.  On $\cm(\cz,\Gamma_Z)$:
\begin{eqnarray*}
\ev_i^*\oh(D)&=&\Le{i};
\end{eqnarray*}
and on $\cm(\ca,\Gamma_A)$:
\begin{eqnarray*}
\Top\otimes\Bot&=&\Split\\
\Top\otimes\ev_i^*L^\vee&=&\Lnt{i}\\
\Bot\otimes\ev_i^*L&=&\Lnb{i}.
\end{eqnarray*}

In section 6, we prove degeneration formulae involving these
line-bundles.  The degeneration formulae are the precise
mathematical statements of the geometric interpretations described
above.

In section 7, we prove a technical result called the {\em Trivial
Cylinder Theorem} which is necessary to the formalism of
\cite{Ka2}.

This paper draws most directly on the Relative Gromov-Witten
Invariants constructed by J. Li \cite{Li1,Li2}. Other approaches
to relative invariants include those of Gathmann \cite{Gath1},
Ionel and Parker \cite{IP}, and A.-M. Li and Ruan \cite{LiRuan}.

Because of the technical nature of this paper, it is not to be
read independently of \cite{Ka2}.  It is our sincere hope that
making this paper available will be useful to other researchers.

I would like to acknowledge the following for valuable
conversations: Y. Eliashberg, A. Gathmann, D. Hain, J. Li, and R.
Vakil. This paper, together with \cite{Ka2} is a revised version
of the author's Ph.D. thesis written under the direction of Y.
Elisahberg.


    \section{Definitions}

\parskip=12pt




\subsection{Schemes with the action of a group}

We need to work in the category of schemes with the action of a
group, or swags.

\begin{definition}\label{swag} {A scheme with the action of a group
is a triple $(X,G,m)$ where $X$ is a scheme, $G$ is a group
scheme, and $m$ is a morphism
$$m:G\times X \srarr X$$
that gives a group action of $G$ on $X$}
\end{definition}

\begin{definition}\label{mswag} A morphism of swags
$(X,G,m)\srarr(X',G',m')$ is a pair $(f,f_*)$ where $f:X\srarr X'$
is a morphism, $f_*:G\srarr G'$ is a group scheme homomorphism and
the following diagram commutes.
$$\xymatrix{
G\times X \ar[r]^<<<<<m\ar[d]_{f_* \times f} & X\ar[d]^f\\
G'\times X' \ar[r]_<<<<<{m'}&X'}$$
\end{definition}

\begin{definition} Given a morphism of swags
$(f,f_*):(X,G,m)\srarr(X',G',m')$, the {\em swag-structure on $X'$
induced by $(f,f_*)$} is the swag $(X',G,m'')$ whose group action
is given by the following composition
$$\xymatrix{
G\times X'\ar[r]^{f_*\times\id} &G'\times X'\ar[r]^>>>>>m &X'}$$
\end{definition}

\begin{definition}\label{vbswag} {A {\em vector bundle} on a swag
$(X,G,m)$ is an equivariant bundle $E$ on $X$}
\end{definition}

\begin{definition}\label{pbswag} {Given a morphism of swags
$$(f,f_*):(X,G,m)\srarr(X',G',m')$$
and a vector bundle $E$, $X'$, {\em the pullback, $f^*E$} is
defined to be the equivariant pullback of $E$ under the
$G$-equivariant map $f:(X,G,m)\srarr (X',G,m'')$ where
$(X',G,m'')$ is the swag-structure on $X'$ induced by $(f,f_*)$.}
\end{definition}

\subsection{Convention for Projectivizations}

If $E$ is a vector-bundle on a scheme $E$, let $\proj_X(E)$ denote
the projectization of $E$ where we use the old-fashioned geometric
notation where points in the projectivization represent lines.  In
this case, if $1_X$ is the trivial bundle on $X$, $\proj_X(E\oplus
1_X)$ is the projective completion of $E$.  The scheme
$$\proj_X(E\oplus 1_X)\setminus E$$
is called the infinity section.

\noindent



    \section{Construction of moduli stacks}



In this section, we construct moduli stacks $\cm(\cz,\Gamma_Z)$,
$\cm(\ca,\Gamma_A)$.  The material in this chapter is a
straightforward adaptation of \cite{Li1} and \cite{Li2}.  While J.
Li does not construct $\cm(\ca,\Gamma_A)$, our construction
directly parallels his.  We do change some notation from
\cite{Li1} to suit our purposes.

Consider a projective manifold $Z$ with a smooth divisor $D$.
$\cm(\cz,\Gamma)$ is the moduli stack of stable maps of curves to
$Z$ relative to $D$.  That is, we look at stable maps where we
specify multiplicities at $D$ together with the usual genus and
degree information.  This data is summarized in a relative graph
$\Gamma$.  Because the target $Z$ may degenerate, we are forced to
introduce $\aff^n$-schemes $Z[n]$ which measure degenerations. The
schemes $Z[n]$ are constructed by an inductive procedure where
$Z[0]=Z$, and the inductive step is similar to deformation to the
normal cone. Now, $Z[n]$ admits a natural $G[n]=(\cs)^n$-action
where we will consider two degenerations of $Z$ equal if they are
related by the $G[n]$-action.  We introduce a stack, $\cz^\rel$
that models families that are locally isomorphic to $Z[n]\srarr
\aff^n$. Following that, we introduce maps of families of curves
to families in $\cz^\rel$.  We impose the conditions of {\em
pre-deformability} and {\em stability} on these maps to ensure
they form a proper Deligne-Mumford stack $\cm(\cz,\Gamma)$. We
recall the requisite results from \cite{Li2} to construct a
tangent obstruction complex and a virtual cycle.

For $\cm(\ca,\Gamma)$, we consider a related geometric situation.
Let $X$ be a projective manifold and $L$ a line bundle on $X$. Let
$P=\proj_X(L\oplus 1_X)$, the projective completion of $L$ which
is a $\proj^1$-bundle. $P$ has two important divisors, $D_0$ and
$D_\infty$, the zero and infinity sections.  We study stable maps
to $P$ relative to $D_0$ and $D_\infty$ where we mod out by a
$\cs$-factor that dilates the fibers.  Again, the target $P$ may
degenerate, we need to introduce a sequence of $\aff^n$-schemes
$A[n]$ where $A[0]=P$. $A[n]$ admits a natural
$\Gt{n}=(\cs)^{n+1}$-action where one of $\cs$ factor comes from
dilating the fibers of $P\srarr X$, and the others are analogous
to those on $Z[n]$.  We introduce a stack of degenerations,
$\ca^\rel$ and a proper Deligne-Mumford stack $\cm(\ca,\Gamma)$ of
maps to $P$ modulo dilation of the fibers.  This stack is called
the stack of {\em maps to rubber} and carries a virtual cycle.

We should explain our top/bottom convention. In $Z$, moving
towards $D$ is considered moving towards the top. In $P$, $D_0$ is
considered the top while $D_\infty$ is the bottom. This slightly
odd convention makes sense in that the most natural choice for
$(X,L)$ is $(D,N_{D/Z})$.  In this case, the zero section of $P$
is identified with $D$ and the normal bundle to $D_0$ in $P$ is
equal to the normal bundle to $D$ in $Z$. Therefore, $D_0\subset
P$ like $D\subset Z$ is on top.

\subsection{The spaces $Z[n]$}

We will study maps to a projective manifold $Z$ relative to a
divisor $D$.  Because of the nature of relative invariants, we
must allow the target $Z$ to degenerate. Let $L=N_{D/Z}$ be the
normal bundle to $D$ in $Z$.  Let $P=\proj(L\oplus 1_D)$ be the
projective completion of $L$.  $P$ has two distinguished divisors,
$D_0$ and $D_\infty$, the zero and infinity sections of $L$.

\begin{definition} Let $_k Z$ be the union of $Z$ with $k$ copies
of $P$, $Z\sqcup_D P_1 \sqcup_D \dots \sqcup_D P_k$, the scheme
given by identifying $D\subset Z$ with $D_\infty \subset P_1$ and
$D_0\subset P_i$ with $D_\infty \subset P_{i+1}$ for
$i=0,1,\dots,k-1$.  $_k Z$ has a distinguished divisor
$D=D_0\subset P_k$.
\end{definition}

We will construct a sequence of pairs of projective manifolds and
smooth divisors, $(Z[n],D[n])$ that map to $\aff^n$ so that the
fiber over different closed points in $\aff^n$ is $(\ _k Z,D)$ for
$k=0,\dots,n$.

\begin{definition} \label{Z0} Let $Z[0]=Z$ and $D[0]=D$.

$(Z[n],D[n])$ is defined inductively.  Let
\begin{enumerate}
\item[(1)] $Z[n]=\Bl_{D[n-1]\times\{0\}}(Z[n-1]\times\aff^1).$
\item[(2)] $D[n]$ is the proper transform of $D[n-1]\times\aff^1$.
\end{enumerate}
\end{definition}

 \begin{definition} \label{nn} {Let
$[n]$ denote the set of integers $\{1,2,\dots,n\}$}.
\end{definition}

\begin{definition} \label{Gn} {Let $G[n]$ denote the group
$(\com^*)^n$.  An element $\sigma\in G[n]$ can be written as an
n-tuple, $\sigma=(\sig{1},\sig{2},\dots,\sig{n})$}
\end{definition}

$G[n]$ acts on $\aff^n$, $G[n]\times\aff^n \srarr \aff^n$ by
$$(\sig{1},\sig{2},\dots,\sig{n})\cdot(t_1,t_2,\dots,t_n) \mapsto
(\sig{1}^{-1}t_1,\sig{2}^{-1}t_2,\dots,\sig{n}^{-1}t_n).$$
$(Z[n],D[n])$ has the following properties:
\begin{enumerate}
\item[(1)] There exists a morphism $p:Z[n]\srarr \aff^n$.

\item[(2)] $Z[n]$ admits a $G[n]$-action fixing $D[n]$ point-wise,
so that $p$ is equivariant.

\item[(3)] $D[n]$ is given as the zero-scheme of a section $s[n]$
of a line-bundle $V[n]$.

\item[(4)] For any order-preserving inclusion
$i:[k]\hookrightarrow[n]$, there is an injective homomorphism
$G[k]\hookrightarrow G[n]$, an inclusions of swags
$\aff^k\hookrightarrow\aff^n$, $Z[k]\hookrightarrow Z[n]$ such
that the following diagram commutes:
$$\xymatrix{
Z[k]\ar[r]\ar[d]_p &Z[m]\ar[d]^p \\
\aff^k\ar[r]&\aff^n}$$

The maps induce an isomorphism $Z[k]\cong
Z[n]\times_{\aff^n}\aff^k$.

\item[(5)] There is a morphism $c:Z[n]\srarr Z$ such that for any
inclusion as in (3), the following diagram commutes
$$\xymatrix{
Z[k]\ar[rr]\ar[dr]_c&&Z[n]\ar[dl]^c&\\
&Z&}$$
\end{enumerate}

Let us address (3). Since $Z$ is smooth $D$ is represented by a
Cartier divisor $(V, s)$ where $V$ is a line-bundle and $s$ is a
section of $V$ whose zero scheme is $D$. $Z[n]=\Proj_{Z[n-1]\times
\aff^1}\left(\bigoplus_{n\geq 0} \ci^n\right)$ for a sheaf of
ideals $\ci$.  Let $\pi:Z[n]\srarr Z[n-1]\times \aff^1\srarr
Z[n-1]$ be the projection. $V[n]=\pi^*V[n-1]\otimes \oh(1)$ and
$s[n]$ is $\pi^*s[n-1]$ considered as a degree one element.
Therefore, we have,

\begin{lemma}\label{pddn} $D[n]$ is represented by a Cartier
divisor of a section $s[n]$ of a line-bundle $V[n]$.
\end{lemma}

The contraction map $c:Z[n]\srarr Z$ is given by the composition
$$c:Z[n]\srarr Z[n-1]\srarr \dots \srarr Z[1]\srarr Z[0].$$

The projection $p:Z[n]\srarr \aff^n$ is defined inductively by the
composition
$$p:Z[n]\srarr Z[n-1]\times\aff^1 \srarr \aff^{n-1}\times\aff.$$

\subsection{Effective Maps}

Let us give a geometric description of $Z[n]$.  The construction
of $Z[n]$ is similar to that of the deformation of the normal cone
in \cite{Fulton}.  Let $N$ be the normal bundle to $D$ in $Z$. Let
$P=\proj_X(N\oplus 1)$ be the projective completion of the total
space of $N$. $Z[1]$ is an $\aff^1$-scheme.  The fibers of this
map over $t\in\aff^1$, $t\neq 0$ is $Z$ while the fiber over $t=0$
is $Z\sqcup_D P$, that is, $Z$ union $P$ along $D\subset Z$ which
is isomorphic to $D_\infty$, the infinity divisor in $P$. The
$G[1]$ group action permutes the fibers over $\com\setminus 0$ and
dilates $P$, considered as a $\proj^1$-bundle.

The fiber of $Z[n]$ over $t\in \aff^n$ is $\ _k Z$.  $k$ is the
number of $0$'s among the coordinates of $t$.  If the $i$th zero
in $t$'s coordinate occurs in the $j$th place, then the $j$th
factor of $\com^*$ in $G[n]$ dilates the fibers of $P_k$ for all
$k\geq i$. We will show this in \ref{degenfibers}

\begin{definition} {Let $H_l\subset\aff^n$ be the hyperplane in
$\aff^n$ given by the equation $t_l=1$}
\end{definition}

Let $i$ be the morphism of swags given by
$$\begin{array}{ccccc}
i&:&\aff^{n-1}=H_l&\hookrightarrow&\aff^n\\
i_*&:&G[n-1]&\srarr & G[n]\\
i_*&:&(\sig{1},\dots,\sig{n-1})&\mapsto&
(\sig{1},\dots,\sig{l-1},1,\sig{l+1},\dots,\sig{n-1})
\end{array}$$

Note that this inclusion makes $Z[n]$ a $G[n-1]$-scheme.  We call
this the {\em standard inclusion}.

\begin{lemma} {$Z[n]\times_{\aff^n} H_l=Z[n-1]$ as
$G[n-1]$-schemes}
\end{lemma}

Let $i:[k]\srarr [n]$ be an order preserving inclusion.  This
induces an inclusion,
$$\aff^k\hookrightarrow\aff^n$$
where is $(s_1,\dots,s_k)$ are the coordinates on $\aff^k$ and
$(t_1,\dots,t_n)$ are the coordinates on $\aff^n$, $t_{i(m)}=s_m$
for $1\leq m\leq k$ and the other coordinates on $\aff^n$ are set
to $1$.  This inclusion can be factored as the inclusion of
hyper-planes.  Given an order-preserving inclusion $i:[k]\srarr
[n]$, there is an inclusion $Z[k]\hookrightarrow Z[n]$, given as
the composition of inclusions as above.

\begin{lemma} Given an order-preserving inclusion $i:[k] \srarr
[n]$, the induced maps fit into the following commutative diagram
$$\xymatrix{
Z[k]\ar[d]_p\ar[r]&Z[n]\ar[d]^p \\
\aff^k\ar[r]&\aff^n}$$
\end{lemma}

We have an isomorphism
$$Z[n]\times_{\aff^n}\aff^k\cong Z[k].$$

The spaces $Z[n]$ are models of families of degenerations of $Z$.

\begin{definition} {If $S$ is any scheme then an {\em effective
family of $(Z,D)$} over $S$ is the triple
$(\tau,\tilde{\cz},\tilde{D})$ given by a morphism
$$\tau:S\srarr\aff^n$$
where
$$\tilde{\cz}=\tau^*Z[n]=Z[n]\times_{\aff^n} S$$
$$\tilde{D}=\tau^*D[n]=D[n]\times_{\aff^n} S$$}
\end{definition}

\begin{definition} Given $\rho:S\srarr G[n]$ and an effective
family $(\tau,\tilde{\cz},\tilde{D})$ over $S$, the {\em action of
$\rho$ on $(\tau,\tilde{\cz},\tilde{D})$} is the family
$(\tau^\rho,\tilde{\cz^\rho},\tilde{D^\rho})$ induced by the map
$$\xymatrix{
\tau^\rho:S\ar[r]^>>>>>{\rho\times\tau}&G[n]\times\aff^n\ar[r]^<<<<<m&\aff^n}$$
\end{definition}

\begin{definition} {Given two effective families,
$$\xi_1=(\tau_1,\tilde{\cz_1},\tilde{D_1})$$
$$\xi_2=(\tau_2,\tilde{\cz_2},\tilde{D_2})$$
associated to morphisms
$$\tau_1:S\srarr \aff^{n_1}$$
$$\tau_2:S\srarr \aff^{n_2}$$
where $n_1\leq n_2$, then an {\em effective map $a$ from $\xi_1$
to $\xi_2$}, $$a:\xi_1\srarr \xi_2$$ consists of
\begin{enumerate}
\item[(1)] An inclusion $i:\aff^{n_1}\srarr \aff^{n_2}$ associated
to an order preserving inclusion $[n_1]\srarr [n_2]$ \item[(2)] A
morphism $\rho:S\srarr G[n_2]$
\end{enumerate}
so that
$$(i\circ \tau_1)^\rho=\tau_2.$$}
\end{definition}

It is clear that the composition of two effective maps is an
effective map.

\begin{definition} {Two effective families,
$$\xi_1=(\tau_1,\tilde{\cz_1},\tilde{D_1})$$
$$\xi_2=(\tau_2,\tilde{\cz_2},\tilde{D_2})$$
are said to be related by an {\em effective arrow}
$$F:\xi_1\srarr \xi_2$$ if there is a {\em dominating} family,
$$\xi=(\tau,\tilde{\cz},\tilde{D})$$
together with effective maps
$$a_1:\xi_1\srarr\xi$$
$$a_2:\xi_2\srarr\xi$$}
\end{definition}

One can easily show that the composition of effective arrows is an
effective arrow by constructing a family which dominates all the
families in the composition.

\begin{lemma} An effective family $(\tau,\tilde{\cz},\tilde{D})$
has a {\em collapsing map} to $(Z,D)$,
$$c_{\tilde{\cz}}:\tilde{\cz}\srarr Z$$
so that $c(\tilde{D})\subset D$
\end{lemma}

\begin{proof} The collapsing is induced from the composition
$$c\circ\tau:S\srarr Z[n]\srarr Z$$
\end{proof}

\begin{lemma} \cite{Li1}, Lemma 4.3 {Let $\xi_i=(\tau_i,\tilde{\cz_i},\tilde{D_i})$ be
effective families over $S$.  Suppose there is an isomorphism
$f:(\tilde{\cz_1},\tilde{D_1})\srarr (\tilde{\cz_2},\tilde{D_2})$
that fits into a commutative diagram
$$\xymatrix{
\tilde{\cz_1}\ar[rr]\ar[dr]_{c\times\pi} & &\tilde{\cz_2}\ar[dl]^{c\times\pi} \\
&Z\times S&}$$
Let $p\in S$.  Then there is an open neighborhood
$U$ of $p$ in $S$ so that over $U$,
$$f|_U:\tilde{\cz_1}\times_S U\srarr \tilde{\cz_2}\times_S U$$
is induced by an effective arrow between $\xi_1\times_S U$ and
$\xi_2\times_S U$.}
\end{lemma}

\begin{definition} \label{forp} If S is a scheme, {\em a family of
relative pairs over $S$} is a triple
$(\tilde{\cz},\tilde{D},\tilde{\pi})$ where $\tilde{\pi}:\cz\srarr
Z\times S$ is a morphism and $\tilde{D}$ is a Cartier divisor of
$\tilde{\cz}$ such that the following property holds: There exists
an open covering $\{U_\alpha\}$ of $S$ such that for all $\alpha$,
$(\tilde{\cz}\times_S U_\alpha,\tilde{D}\times_S U_\alpha)$ is
isomorphic to an effective family
$(\tau_\alpha,\tilde{\cz_\alpha},\tilde{D_\alpha})$ over
$U_\alpha$ by a morphism $f$ such that the following diagram
commutes
$$\xymatrix{
(\tilde{\cz}\times_S U_\alpha,\tilde{D}\times_S
U_\alpha)\ar[rr]\ar[dr]_{\tilde{\pi}} && (\tilde{\cz_\alpha},\tilde{D_\alpha})\ar[dl]^{c\times\pi} \\
& Z\times S &}$$
\end{definition}

\begin{definition} Let $\cz^\rel$ be the category whose objects
are relative families of $(Z,D)$ over schemes.  For
$\xi_1,\xi_2\in\Obj(\cz^\rel)$ such that $\xi_i$ is a family over
$S_i$, $\Mor(\xi_1,\xi_2)$ consists of pairs $(h,f)$ where
$$h:S_1\srarr S_2$$
and $f$ is an isomorphism fitting into the following commutative
diagram
$$\xymatrix{
\xi_1\ar[rr]^f\ar[dr]&&h^*\xi_2\ar[dl]\\
&Z\times S_1&}$$

Define $\cp:\cz^\rel\srarr\Sch$ to be the functor that sends
families over $S$ to $S$.  $(\cz^\rel,\cp)$ is a groupoid.
\end{definition}

\begin{proposition} {$(\cz^\rel,\cp)$ is a stack.}
\end{proposition}


\subsection{The spaces $A[n]$}

We construct a sequence of spaces $A[n]$, called {\em the rubber
target spaces} that are used to study curves in the
projectivization of a line bundle $L$ over a projective manifold
$X$. These spaces will be used to construct the rubber invariants.
These spaces have properties analogous to the $Z[n]$'s.

Let $X$ be a projective manifold and $L$ a line bundle on $X$. Let
$P=\proj_X(L\oplus 1_X)$, and let $X_0$ and $X_\infty$ denote the
zero and infinity sections.  We study stable maps to $P$ relative
to $X_0$ and $X_\infty$ where we mod out by a $\cs$-factor that
dilates the fibers.  Again, the target $P$ may degenerate.

\begin{definition} Let $_k P$ be the union of $k+1$ copies of $P$,
$$_k P =P_0 \sqcup_X P_1 \sqcup_X \dots \sqcup_X P_k$$
gluing $X_0\subset P_i$ to $X_\infty\subset P_{i+1}$ for
$i=0,\dots,k-1$.
\end{definition}

$_k P$ has distinguished divisors $D_\infty=X_\infty\subset P_0$
and $D_0=X_0 \subset P_k$.  We will construct a sequence of
projective manifolds and divisors $(A[n],D_0[n],D_\infty[n])$.
$A[n]$ will map to $\aff^n$ with fiber over closed points equal to
$_k P$ for varying $k$.

\begin{definition} \label{Gtn} {Let $\Gt{n}$ denote the group
$(\com^*)^{n+1}$.  An element $\sigma\in \Gt{n}$ can be written as
an $n+1$-tuple, $\sigma=(\sig{0},\sig{1},\dots,\sig{n})$}
\end{definition}

In what follows, we will have $\Gt{n}$ act on $\aff^n$,
$\Gt{n}\times\aff^n \srarr \aff^n$ as
$$(\sig{0},\sig{1},\sig{2},\dots,\sig{n})\cdot(t_1,t_2,\dots,t_n)
\mapsto (\sig{1}^{-1}t_1,\sig{2}^{-1}t_2,\dots,\sig{n}^{-1}t_n).$$

\begin{definition} The spaces $(A[n],D_0[n],D_\infty[n])$ are
defined as $\Gt{n}$-schemes as follows
\begin{enumerate}
\item[(1)] $A[0]=P=\proj_X(L\oplus 1_X)$ where

\subitem(a) $D_0[0]$ is the zero-section in $A[0]$, which is seen
as the projective completion of $L$, alternatively as a divisor of
the line bundle $\oh_P(1)\otimes \pi^*L$.

\subitem(b) $D_\infty[0]$ is the infinity section of $A[0]$, that
is the divisor of the line bundle $\oh_P(1)$

\subitem(c) The $\Gt{0}=\com^*$-action is given by
$$\sig{0}\cdot [l:t] = [\sig{0}l:t]$$

\item[(2)] $A[n]$ is defined inductively from $A[n-1]$ as follows.
Consider $A[n-1]\times\aff^1$ under the group action
$$(\sig{0},\sig{1},\dots,\sig{n})\cdot(z,t_n)=((\sig{0},\sig{1},\dots,\sig{n-1})\cdot
z,\sig{n}^{-1}t_n)$$

Let $A[n]=\Bl_{D_0[n-1]\times\{0\}}(A[n-1]\times\aff^1)$.

\subitem(a) $D_0[n]$ is the proper transform of
$D_0[n-1]\times\aff^1$.

\subitem(b) $D_\infty[n]$ is the proper transform of
$D_\infty[n-1]\times\aff^1$.

\end{enumerate}
\end{definition}

$(A[n],D_0[n],D_\infty[n])$ has the following properties
\begin{enumerate}
\item[(1)] There exists a morphism $p:A[n]\srarr \aff^n$.

\item[(2)] $A[n]$ admits a $\Gt{n}$-action fixing $D_0[n]$ and
$D_\infty[n]$ making $p$ equivariant.

\item[(3)] For any order-preserving inclusion
$[k]\hookrightarrow[n]$, there is an injective homomorphism
$\Gt{k}\hookrightarrow \Gt{n}$, an inclusions of swags
$\aff^k\hookrightarrow\aff^n$, $A[k]\hookrightarrow A[n]$,  such
that the following diagram commutes:
$$\xymatrix{
A[k]\ar[r]\ar[d]_p & A[n]\ar[d]^p \\
\aff^k\ar[r] & \aff^n}$$
 Moreover, the maps induce an isomorphism
$A[k]\cong A[n]\times_{\aff^n}\aff^k$.

\item[(4)] There is a natural projection
$$\pi:A[n]\srarr X$$
so that for any inclusion as in (3), the following diagram
commutes
$$\xymatrix{
A[k]\ar[rr]\ar[dr]_{\pi} & & A[n] \ar[dl]^{\pi}\\
& X &}$$

\item[(5)] There is a morphism of swags $t:A[n]\srarr A[0]$ where
$t_*:\Gt{n}\srarr\Gt{0}$ is given by
$$t_*:(\sig{0},\sig{1},\dots,\sig{n})\srarr
(\sig{0}\sig{1}\dots\sig{n}).$$ For any inclusion as in (3), the
following diagram commutes
$$\xymatrix{
A[k] \ar[rr]\ar[dr]_t&& A[n]\ar[dl]^t \\
& A[0] &}$$ $t$ is called the {\em top} morphism.

\item[(6)] There is a morphism of swags $b:A[n]\srarr A[0]$ where
$b_*:\Gt{n}\srarr\Gt{0}$ is given by
$$b_*:(\sig{0},\sig{1},\dots,\sig{n})\srarr (\sig{0}).$$
For any inclusion as in (3), the following diagram commutes
$$\xymatrix{
A[k] \ar[rr]\ar[dr]_b&& A[n]\ar[dl]^b \\
& A[0] &}$$
 $b$ is called the {\em bottom} morphism.

\end{enumerate}

Properties (1-4) follow from proofs analogous to those for $Z[n]$.

\begin{definition} {The map $b:A[n]\srarr A[0]$ is given as the
composition of blow-downs and projections
$$b:A[n]\srarr A[n-1]\times \aff^1 \srarr A[n-1] \srarr \dots
\srarr A[0] \times \aff^1 \srarr A[0].$$}
\end{definition}

\begin{definition} \label{topmap} We can construct $A[n]$ by repeatedly blowing up $D_\infty$
instead of $D_0$.  The map $t:A[n]\srarr A[0]$ is given as a
composition of blow-downs and projections analogous to $b$.
\end{definition}

The terminology for the top and bottom maps is as follows.  Given
a point $x \in \aff^n$, the fiber in $A[n]$ over $x$ is
$$_l P=P_0\sqcup_X P_1 \sqcup_X\dots\sqcup_X P_l$$
where $l$ is the number of zeroes among $x$'s coordinates.  The
map $b$ restricted to the fiber over $x$ is the map $b:\ _l
P\srarr P$ so that
\begin{enumerate}
\item[(1)] On $P_0$, it is the identity $P_0 \srarr P$.

\item[(2)] On $P_i$ for $i\geq 1$, it is the projection $P_i\srarr
X=D_0\subset P$.
\end{enumerate}
Likewise, the map $t:A[n]\srarr A[0]$ restricts to $_l P$ as
\begin{enumerate}
\item[(1)] On $P_l$, it is the identity $P_l \srarr P$.

\item[(2)] On $P_i$ for $i\leq l-1$, it is the projection
$P_i\srarr X=D_\infty\subset P$.
\end{enumerate}

The definition for effective arrows for $A[n]$ is analogous to
that of $Z[n]$ with $Z$'s replaced replaced by $A$'s and $G[n]$'s
replaced by $\Gt{n}$'s. We have {\em families of rubber} modelled
on maps to $S\srarr \aff^n$ analogous families of relative pairs
\ref{forp}

In particular, we have

\begin{lemma} {Let
$\xi_i=(\tau_i,\tilde{\ca_i},\widetilde{D_0}_i,\widetilde{D_\infty}_i)$
be rubber families over $S$.  Suppose there is an isomorphism
$f:(\tilde{\ca_1},\widetilde{D_0}_1,\widetilde{D_\infty}_1)\srarr
(\tilde{\ca_2},\widetilde{D_0}_2,\widetilde{D_\infty}_2)$ that
fits into a commutative diagram
$$\xymatrix{
\tilde{\ca_1}\ar[rr]\ar[dr]_{\pi\times q} & & \tilde{\ca_2}\ar[dl]^{\pi\times q} \\
& X\times S &}$$
Let $p\in S$.  Then there is an open neighborhood
$U$ of $p$ in $S$ so that over $U$,
$$f|_U:\tilde{\ca_1}\times_S U\srarr \tilde{\ca_2}\times_S U$$
is induced by an effective arrow between $\xi_1\times_S U$ and
$\xi_2\times_S U$.}
\end{lemma}

\begin{definition} \label{carel} Let $\ca^\rel$ be the category
whose objects are relative families of $(X,L)$, modelled on
$A[n]\srarr \aff^n$ over schemes. For
$\xi_1,\xi_2\in\Obj(\ca^\rel)$ such that $\xi_i$ is a family over
$S_i$, $\Mor(\xi_1,\xi_2)$ consists of pairs $(h,f)$ where
$$h:S_1\srarr S_2$$
and $f$ is an isomorphism fitting into the following commutative
diagram

$$\xymatrix{
\xi_1\ar[rr]^f\ar[dr]&&h^*\xi_2\ar[dl]\\
&X\times S_1&}$$

Define $\cp:\ca^\rel\srarr\Sch$ to be the functor that sends
families over $S$ to $S$.  $(\ca^\rel,\cp)$ is a groupoid.
\end{definition}

\begin{proposition} {$(\ca^\rel,\cp)$ is a stack.}
\end{proposition}

\subsection{Degenerate Fibers}\label{degenfibers}

We identify the fibers of $Z[n]\srarr\aff^n$.  Let $N=N_{D/Z}$ be
the normal bundle to $D$ in $Z$.  Let $A[n]$ be the spaces
constructed from $(D,N)$.

Let $K_k\subset \aff^n$ denote the subset of $\aff^n$ given by the
equation $t_k=0$.  We need to understand the fiber of $K_k$ under
$p:Z[n]\srarr \aff^n$.

Consider $Z[k-1]\times\aff^{n-k}$ with the $G[n]$- action given by
$$(\sig{1},\dots,\sig{n})\cdot (z,(t_{k+1},\dots,t_n))\mapsto
((\sig{1},\dots,\sig{k-1})\cdot
z,(\sig{k+1}^{-1}t_{k+1},\dots,\sig{n}^{-1}t_n)$$  Note that this
scheme has a divisor $D[k-1]\times\aff^{n-k}$.

Let us also consider $\aff^{k-1}\times A[n-k]$ with the
$G[n]$-action given by
$$(\sig{1},\dots,\sig{n})\cdot (t_1,\dots,t_{k-1},z)\mapsto
((\sig{1}^{-1}t_1,\dots,\sig{k-1}^{-1}t_{k-1}),(\sig{1}\dots\sig{k},\sig{k+1},\dots,\sig{n})\cdot
z)$$  Note that this scheme has a divisor given by
$\aff_{k-1}\times D_\infty[n-k]$.  Note that the two divisors are
isomorphic to $D\times \aff^{n-1}$ and that they have the same
$G[n]$-action.

\begin{lemma} {The fiber over $K_k$ is given by the following
isomorphism of $G[n]$-schemes.
$$Z[n]\times_{\aff^n} K_k = (Z[k-1]\times \aff^{n-k})
\sqcup_{D\times\aff^{n-1}} (\aff^{k-1}\times A[n-k])$$ where
$\sqcup_D$ denotes union identifying the divisor on each scheme.}
\end{lemma}

\begin{proof} This follows from standard facts about blow-ups and by
induction. \end{proof}

By induction, we can give a description of fibers over closed
points.  Let $t\in \aff^n$ have $l$ zeroes among its coordinates.
Then the fiber of $Z[n]$ over $t$ is $_l Z$.

An analogous result holds for $A[n]$. Let us put the following
$\Gt{n}$-structure on $A[k-1]\times\aff^{n-k}$
$$(\sig{0},\sig{1},\dots,\sig{n})\cdot (z,(t_{k+1},\dots,t_n))\mapsto
((\sig{0},\sig{1},\dots,\sig{k-1})\cdot
z,(\sig{k+1}^{-1}t_{k+1},\dots,\sig{n}^{-1}t_n)$$  Note that this
scheme has a divisor $D_0[k-1]\times\aff^{n-k}$. Let us also
consider $\aff^{k-1}\times A[n-k]$ with the $\Gt{n}$-action given
by
$$(\sig{0},\sig{1},\dots,\sig{n})\cdot (t_1,\dots,t_{k-1},z)\mapsto
((\sig{1}^{-1}t_1,\dots,\sig{k-1}^{-1}t_{k-1}),(\sig{0}\sig{1}\dots\sig{k},\sig{k+1},\dots,\sig{n})\cdot
z)$$  Note that this scheme has a divisor given by
$\aff^{k-1}\times D_\infty[n-k]$.

\begin{lemma} {The fiber over $K_k$ is given by the following
isomorphism of $\Gt{n}$-schemes.
$$A[n]\times_{\aff^n} K_k = (A[k-1]\times \aff^{n-k}) \sqcup_D
(\aff^{k-1}\times A[n-k])$$}
\end{lemma}

Let $t\in \aff^n$ have $l$ zeroes among its coordinates.  Then the
fiber of $A[n]$ over $t$ is $_l P$.

\subsection{Pre-deformable Families}

Note that $\Sing(_k Z)$, the singular locus of $_k Z$ is the
disjoint union of $k$ copies of $D$, which we label
$D_1,\dots,D_{k-1}$ where $D_i=D_\infty\subset P_i$.

\begin{definition} A morphism $f:C\srarr\ _k Z$ is said to be {\em
pre-deformable} if $f^{-1}(D_i)$ is the union of nodes so that for
$p\in f^{-1}(D_i)$ ($i=1,2,\dots,k)$, the two branches of the node
map to different irreducible component of $_k Z$ and that the
order of contact to $D_i$ are equal.
\end{definition}

An obvious analogous definition exists for morphisms to $\ _k P$.

There is a precise notion of pre-deformable families of morphisms
given in terms of local models.  See \cite{Li1} for details.

\begin{definition} A {\em morphism of a pre-stable family with
marked points} over a scheme $S$ consists
\begin{enumerate}
\item[(1)] A family of targets, $$Y\in\Obj(\cz^{rel}(S))$$

\item[(2)] A pre-stable family $\ex\srarr S$.

\item[(3)] A morphism $f$ fitting into
$$\xymatrix{
\ex\ar[rr]^f \ar[dr]& &Y\ar[dl]\\
&S&}$$

\item[(4)] Morphisms $\gamma_1,\dots,\gamma_k:S\srarr \ex$
\end{enumerate}
so that in a neighborhood $S_\alpha\srarr S$ where $\cz\times_S
S_\alpha$ is given by a morphism to $\aff^n$, the diagram looks as
follows
$$\xymatrix{
\ex \ar[d] \ar[r]^f & Z[n]\ar[d]\\
S_\alpha \ar@<1ex>[u]^{\gamma_1,\dots,\gamma_k}\ar[r] & \aff^n}$$

where

\begin{enumerate}
\item[(1)] $f$ is a pre-deformable morphism.

\item[(2)] $\gamma_1,\dots,\gamma_k$ are disjoint section of
$\ex\srarr S$ whose image does not intersect the nodes in the
fibers of $\ex\srarr S$.
\end{enumerate}
\end{definition}

The image of $\gamma$ are marked points.  We will often be
distinguishing a particular marked point and write the following
diagram
$$\xymatrix{
\ex \ar[d] \ar[r]^f & Z[n]\ar[d]\\
S \ar@<1ex>[u]^{\gamma_1}\ar[r] & \aff^n}$$ where
$\gamma_2,\dots,\gamma_k$ are present, but are suppressed in our
notation.  We will often refer to a {\em distinguished marked
point} in the above sense.

\begin{definition} Two morphisms of pre-stable families over $S$,
indexed by $i=1,2$
$$\xymatrix{
\ex_i \ar[d] \ar[r]^{f_i} & \cz_i\ar[dl]\\
S \ar@<1ex>[u]^{\gamma_{i,1},\dots,\gamma_{i,k}} &}$$ are
isomorphic if there exists an isomorphism over $S$
$$\rho:\ex_1\srarr\ex_2$$
and an arrow in $\cz^{rel}(S)$
$$\tau:\cz_1\srarr \cz_2$$
so that the following diagram commutes
$$\xymatrix{
\ex_1\ar[d]_\rho\ar[r]^{f_1}&\cz_1\ar[d]^\tau\\
\ex_2\ar[r]_{f_2}&\cz_2}$$ and
$$\rho\circ\gamma_{1,j}=\gamma_{2,j}$$
for $j=1,\dots,k$
\end{definition}

\begin{definition} A morphisms of families over $S$ is said to be
{\em stable} if for every point $s\in S$, the automorphisms of the
family restricted to $s$ is finite.
\end{definition}

As in absolute Gromov-Witten theory, a contracted curve (one
mapped to a point) occurring as a component in a family is
unstable unless it has genus greater than or equal to two, at
least one marked point if genus is equal to $1$, at least three
marked points if genus is equal to $0$.  In the relative case, we
have to consider the phenomenon of curves mapping to a degenerated
target $\ _l Z$ where we have a $\com^*$ dilating the fiber of
each copy of $P$.  It is shown in (\cite{Li1}, lemma 3.2) that a
morphism $f:C\srarr \ _l Z$ that does not involve any unstable
contracted components is stable if and only if for every $i$,
there is an irreducible component of $C$ mapped into $P_i$ by $f$
that is not a {\em trivial component}, that is, a rational curve
lying in a fiber of $P\srarr X$, totally ramified at points
mapping to $D_0$ and $D_\infty$.

Analogous definitions hold for morphisms to $A[n]$.  For a map
$f:C\srarr \ _l P$, if $f$ does not involve any unstable
contracted components, then $f$ is stable if and only if for every
$i$, some irreducible component of $C$ that is not a trivial
component is mapped into $P_i$.

\subsection{The Moduli Stacks}

We need to specify the appropriate data for the moduli stacks
$\cm(\cz,\Gamma_Z)$, $\cm(\ca,\Gamma_A)$.  Let us begin with
$\cm(\cz,\Gamma)$. Here we consider stable maps to $Z$ with
specified tangency to $D$ together with marked points, called {\em
interior marked points} whose image is not mapped to $D$ by $f$.
We will assume that the points that are mapped to $D$ by $\ex$ are
also marked. Those marked points will be called {\em boundary
marked points}.

\begin{definition} A {\em relative graph} $\Gamma$ is the
following data:
\begin{enumerate}
\item[(1)] A finite set of vertices $V(\Gamma)$

\item[(2)] A genus assignment for each vertex
$$g:V(\Gamma)\srarr \zeeg$$

\item[(3)] A degree assignment for each vertex
$$d:V(\Gamma)\srarr B_1(Z)=A_1(Z)/\tilde{ }_{\text{alg}}$$
that assigns the class of a curve modulo algebraic equivalence to
each vertex.

\item[(4)] A finite set $R=\{1,\dots,r\}$ labelling boundary
marked points together with a function assigning boundary marked
points to vertices
$$a_R:R\srarr V(\Gamma)$$

\item[(5)] A {\em multiplicity assignment} for each boundary
marked point
$$\mu:R\srarr \zee_{\geq 1}$$

\item[(6)] A finite set $M=\{1,\dots,m\}$ labelling interior
marked points together with an assignment to vertices
$$a_M:M\srarr V(\Gamma)$$
\end{enumerate}
\end{definition}

\begin{definition} Given two relative graphs $\Gamma$, $\Gamma'$
are said to be isomorphic if there is a bijection
$$q:V(\Gamma)\srarr V(\Gamma')$$
such that
\begin{enumerate}
\item[(1)] $g(v)=g'(q(v))$ \item[(2)] $d(v)=d'(q(v))$ \item[(3)]
$a_R'(i)=q(a_R(i))$ \item[(4)] $\mu_R'=\mu_R$ \item[(5)]
$a_M'(i)=q(a_M(i))$.
\end{enumerate}
\end{definition}

\begin{definition} Let $\Gamma$ be a relative graph.  A {\em
morphism of families of type $\Gamma$} over $S$ consists of a
stable morphism of a pre-stable family with marked points where
$\cz\in\Obj(\cz^\rel(S))$

$$\xymatrix{
\ex \ar[d] \ar[r]^f & \cz\ar[dl]\\
S \ar@<1ex>[u]^{\gamma_1,\dots,\gamma_m,\delta_1,\dots,\delta_r}
&}$$ such that for any closed point $s\in S$, the fiber
$\ex_s=\ex\times_S s$ obeys
\begin{enumerate}
\item[(1)] $\ex_s$ can be written as a disjoint union of
pre-stable curves $(\ex_s)_v\srarr S$
$$\ex=\coprod_{v\in V(\Gamma)} (\ex_s)_v$$

\item[(2)] $(\ex_s)_v$ is a connected curve of arithmetic genus
$g(v)$.

\item[(3)] The map
$$(c\circ f_s):(\ex_s)_v \srarr \cz \srarr Z$$
has $(c\circ f_s)_*((\ex_s)_v)_s=d(v)$.  Note that the map
$\cz\srarr Z$ is induced from the map $Z[n]\srarr Z$.

\item[(4)] $\gamma_i(s)\subset (\ex_s)_v$ for $v=a_M(i)$.  These
are the interior marked points.  \label{bodymarkedpoint}

\item[(5)] $\delta_i(s)\subset (\ex_s)_v$ for $v=a_R(i)$.  These
are the boundary marked points.  \label{bdrymarkedpoint}

\item[(6)] There is the following multiplicity condition
$$f^*D[n]=\sum_{i\in R} \mu(i)\delta_i(s)$$
\end{enumerate}
\end{definition}

\begin{definition} \label{nicefamily} A morphism of families of
type $\Gamma$ is said to be a {\em nice family} if the target
$\cz$ is an effective family, that is the the morphism be
expressed as
$$\xymatrix{
\ex\ar[r]\ar[r]^f\ar[d]& Z[n]\ar[d]\\
 S \ar@<1ex>[u]^{\gamma_1,\dots,\gamma_m,\delta_1,\dots,\delta_r}\ar[r]
 &\aff^n}$$
\end{definition}

\begin{definition} {The category $\cm(\cz,\Gamma)$ is the groupoid
over (Sch) so that for a scheme $S$, the objects of
$\cm(\cz,\Gamma)(S)$ are morphisms of families of type $\Gamma$.
Given a morphism $\sigma:S\srarr T$, families $\xi_1\in
\cm(\cz,\Gamma)(S)$, $\xi_2\in \cm(\cz,\Gamma)(T)$, an arrow
$\xi_1\srarr\xi_2$ is an isomorphism over $S$ of $\xi_1$ with
$\sigma^*\xi_2$.}
\end{definition}

\begin{theorem} \cite{Li1} $\cm(\cz,\Gamma)$ is a proper Deligne-Mumford
stack.
\end{theorem}

Given a family
$$\xymatrix{
\ex\ar[r]\ar[r]^f\ar[d]& Z[n]\ar[d]\\
 S \ar@<1ex>[u]^{\gamma_1,\dots,\gamma_m,\delta_1,\dots,\delta_r}\ar[r]
 &\aff^n}$$
we have morphisms $S\srarr Z$ given by
$$(c\circ f\circ \gamma_i):S\srarr \ex\srarr Z[n]\srarr Z$$
$$(c\circ f\circ \delta_i):S\srarr \ex\srarr D[n]\srarr D$$
These maps extend to $\cm(\cz,\Gamma)$ giving evaluation maps
$\ev_i:\cm(\cz,\Gamma)\srarr Z$ at the interior and boundary
marked points.

\begin{definition} {The {\em evaluation map} on
$\cm(\cz,\Gamma)$ is
$$\Ev:\cm(\cz,\Gamma)\srarr Z^n\times D^r$$}
\end{definition}

\begin{definition} An \'{e}tale nice family
$S\srarr\cm(\cz,\Gamma_Z)$ is said to be a {\em nice chart}.
\end{definition}

\begin{theorem} (\cite{Li1}, Theorem 3.10) $\cm(\cz,\Gamma_Z)$ has an atlas that is a union of
nice charts.\end{theorem}

The case for $cm(\ca,\Gamma)$ is analogous with only a few
modifications.
\begin{definition} {A {\em rubber graph} $\Gamma$ is the
following data:
\begin{enumerate}
\item[(1)] a finite collection of vertices $V(\Gamma)$

\item[(2)] A genus assignment for each vertex $$g:V(\Gamma)\srarr
\zeeg.$$

\item[(3)] A degree assignment for each vertex
$$d:V(\Gamma)\srarr B_1(X)=A_1(X)/\tilde{ }_{\text{alg}}$$
that assigns the class of a curve modulo algebraic equivalence to
each vertex.

\item[(4)] A finite sets $R_0=\{1,\dots,r_0\}$,
$R_\infty=\{1,\dots,r_\infty\}$ labelling boundary marked points
together with a function assigning boundary marked points to
vertices
\begin{eqnarray*}
a_0&:&R_0\srarr V(\Gamma)\\
a_\infty&:&R_\infty\srarr V(\Gamma)
\end{eqnarray*}

\item[(5)] A {\em multiplicity assignment} for boundary marked
points
\begin{eqnarray*}
\mu^0&:&R_0\srarr \zee_{\geq 1}\\
\mu^\infty&:&R_\infty\srarr \zee_{\geq 1}
\end{eqnarray*}

\item[(6)] A finite set $M=\{1,\dots,m\}$ labelling interior
marked points together with an assignment to vertices
$$a_M:M\srarr V(\Gamma)$$
\end{enumerate}
}
\end{definition}

Definitions of morphisms of rubber families are analogous with
$Z$'s replaced with $A$'s and the following modifications. The
degree assignment is
$$d(v)\in B_1(X)$$
Instead of sections $\delta_1,\dots,\delta_r$, we have sections
\begin{eqnarray*}
\delta^0_1,\dots,\delta^0_{r_0}&:&S\srarr\ex\\
\delta^\infty_1,\dots,\delta^\infty_{r_\infty}&:&S\srarr\ex
\end{eqnarray*}
so that
\begin{eqnarray*}
f^*D_0[n]&=&\sum_{i\in R_0} \mu^0(i)\delta^0_i(s)\\
f^*D_\infty[n]&=&\sum_{i\in R_\infty}
\mu^\infty(i)\delta^\infty_i(s)
\end{eqnarray*}

\begin{theorem} For a rubber graph $\Gamma$,
$\cm(\ca,\Gamma)$ is a proper Deligne-Mumford stack.
\end{theorem}

We have analogous evaluation maps $\ev_i$ at the interior and
boundary marked points (mapping to $D_0$ and $D_\infty$).

\begin{definition} {The {\em evaluation map} on
$\cm(\ca,\Gamma)$ is
$$\Ev:\cm(\ca,\Gamma)\srarr X^n\times X^{r_0}\times X^{r_\infty}$$}
\end{definition}

Now, we introduce some notation which will be essential in the
sequel.
\begin{definition} \label{extendedcomponents} A map $f:C\srarr\ _k
Z$ in $\cm(\cz,\Gamma)$ is said to be {\em split} if $k\geq 1$.
The irreducible components of $C$ that are mapped to $P_i\subset\
_k Z$ are said to be {\em extended components}.
\end{definition}
An analogous situation occurs for $\cm(\ca,\Gamma)$.

\begin{definition} \label{splitmaps} A {\em split map} in $\cm(\ca,\Gamma)$ is
a map $f:C\srarr\ _k P$ where $k\geq 1$, that is, the target is
not smooth.
\end{definition}

\begin{definition} \label{topcomponent}
For a map $f:C\srarr\ _k P$ in $\cm(\ca,\Gamma)$, the irreducible
components of $C$ that are mapped to $P_k$ are said to be the {\em
top components} while the components of $C$ that are mapped to
$P_0$ are said to be the {\em bottom components}.
\end{definition}

\subsection{Gluing Moduli Stacks}

Consider a projective manifold $Z$, together with a smooth divisor
$D$.  We will consider the relative moduli stack
$\cm(\cz,\Gamma_Z)$ corresponding to $(Z,D)$.  Let us look at the
rubber moduli space corresponding to $(X=D,L=N_{D/Z})$ where
$N_{D/Z}$ is the normal bundle to $D$ in $Z$.  $\cm(\ca,\Gamma_A)$
corresponds to certain maps that are added to maps in
$\cm(\cz,\Gamma_Z)$ to compactify.  There is a specific way to
join $\cm(\cz,\Gamma_Z)$ to $\cm(\ca,\Gamma_A)$ if some conditions
are met. Likewise we can join some $\cm(\ca,\Gamma_b)$ to
$\cm(\ca,\Gamma_t)$ if similar conditions are met.  We make these
conditions precise below.

\begin{definition} Let
$\Gamma_Z$ be a relative graph and $\Gamma_A$ be a rubber graph.
Suppose that $L:RZ\srarr RA_\infty$ is a bijection from the
labelling sets for boundary marked points in $\Gamma_Z$ to the
labelling sets for boundary marked points mapping to $D_\infty$ in
$\Gamma_A$ so that
$$\mu_Z(q)=\mu^\infty_A(L(q)).$$
Let
$$J:M_Z\sqcup M_A\srarr \{1,\dots,|M_Z|+|M_A|\}$$
be a bijection between the labelling sets of the interior marked
points and a set of $|M_Z|+|M_A|$ elements.  We call the data
$(\Gamma_A,\Gamma_Z,L,J)$ a {\em graph-join quadruple}.
\end{definition}
Colloquially, we've matched boundary marked points on $\Gamma_Z$
and $\Gamma_A$ with the same multiplicity.

\begin{definition} \label{graphjoin}
Define the {\em graph join} $\Gamma_A*_{L,J}\Gamma_Z$ to be the
following relative graph.  Let the graph $\Delta$ be obtained by
taking as vertices the vertices of $\Gamma_Z$ and $\Gamma_A$ and
for every $q\in RZ$, place an edge between the vertices
corresponding to $q$ and $L(q)$. Let $\Gamma_A*_{L,J}\Gamma_Z$ be
given as follows.  The vertices of
$\Gamma=\Gamma_A*_{L,J}\Gamma_Z$ are the connected components of
$\Delta$. Let $b_Z:V(\Gamma_Z)\srarr V(\Gamma)$, and
$b_A:V(\Gamma_A)\srarr V(\Gamma)$ be the functions taking vertices
of $\Gamma_A$ and $\Gamma_A$ to the components in $\Delta$
containing them. For $v\in V(\Gamma)$, let $\Delta_v$ be the
connected component of $\Delta$ corresponding to $v$.  Define the
data for $\Gamma_A*_{L,J}\Gamma_Z$ as follows:
\begin{enumerate}
\item[(1)] $g(v)=(\sum_{w\in b_Z^{-1}(v)} g(w))+(\sum_{w\in
b_A^{-1}(v)} g(w))+\dim(h^1(\Delta_v))$

\item[(2)] $d(v)=(\sum_{w\in b_Z^{-1}(v)} d(w))+(\sum_{w\in
b_A^{-1}(v)} i_*d(w))$ where $i:X\srarr Z$ is the inclusion and
$i_*:B_1(X)\srarr B_1(Z)$ is the induced map.

\item[(3)] $R=RA_0$ with $a_R:R\srarr V(\Gamma)$ given by
$$a_R=b_A \circ a_0$$

\item[(4)] $\mu:R\srarr \zee_{\geq 1}$ given by
$$\mu_R=\mu^0$$

\item[(5)] $M=\{1,\dots,|M_Z|+|M_A|\}$ with assignment function
$a:M\srarr V(\Gamma)$ given for $k\in J(M_Z)$ by
$$a(k)=b_Z \circ a_{MZ}\circ J^{-1}$$
while for $k\in J(M_A)$ by
$$a(k)=b_A \circ a_{MA}\circ J^{-1}$$
\end{enumerate}
\end{definition}

Given $(\Gamma_Z,\Gamma_A,L,J)$ as above, consider the evaluation
map at the boundary marked points on $\cm(\cz,\Gamma_Z)$ followed
by a map $L_*:D^r\srarr D^r$ which reorders the products of $D^r$
according to $L$:
$$L_*\circ Ev_R:\cm(\cz,\Gamma_Z)\srarr D^r\srarr D^r$$
and the evaluation map at the boundary marked points mapping to
$D_\infty\cong X$ on $\cm(\ca,\Gamma_A)$,
$$\Ev_{R_\infty}:\cm(\ca,\Gamma_A)\srarr D^r.$$

\begin{theorem} \cite{Li1} There is a morphism
$$\Phi_{\Gamma_Z,\Gamma_A,L,J}:\cm(\ca,\Gamma_A)\times_{D^r} \cm(\cz,\Gamma_Z)\srarr \cm(\cz,\Gamma_A*_{L,J}\Gamma_Z).$$
\end{theorem}

\begin{proof}
We look at the morphisms in charts and glue curves together to
form nodes.
\end{proof}

\begin{definition} \label{zuniona} {Define the stack
$\cm(\ca\sqcup\cz,\Gamma_A\sqcup_{L,J} \Gamma_Z)$ as the image
stack of $\Phi$ in $\cm(\cz,\Gamma_A*_{L,J} \Gamma_Z)$.}
\end{definition}

Consider $RZ$ as part of the data of $\Gamma_Z$.  An automorphism
of $RZ$ is a permutation
$$\sigma:RZ\srarr RZ$$
so that $\mu_Z(\sigma(i))=\mu_Z(i)$ and
$a_{RZ}(\sigma(i))=a_{RZ}(i)$.  The group of all such
automorphisms is denoted by $\Aut_{\Gamma_Z}(RZ)$.  Likewise, we
may define $\Aut_{\Gamma_A}(RA_0)$ and
$\Aut_{\Gamma_A}(RA_\infty)$. Given $L:RZ\srarr RA_\infty$, we may
define $\Aut_{\Gamma_A,\Gamma_Z,L}(RZ,RA_\infty)$ as the subgroup
of $\Aut_{\Gamma_Z}(RZ)\times\Aut_{\Gamma_A}(RA_\infty)$ such that
for $(\sigma,\tau)\in
\Aut_{\Gamma_Z}(RZ)\times\Aut_{\Gamma_A}(RA_\infty)$ we have
$L(\sigma(i))=\tau(L(i))$ for $1\leq i\leq |RZ|$.

\begin{lemma} (\cite{Li1}, Proposition 4.13) {$\Phi$ is finite and \'{e}tale onto its image.  It has degree equal
to
$$|\Aut_{\Gamma_A,\Gamma_Z,L}(RZ,RA_\infty)|$$ at every integral
substack of $\cm(\ca\sqcup\cz,\Gamma_A \sqcup_{L,J} \Gamma_Z)$.}
\end{lemma}

\begin{definition} \label{joineqvlnt} {Two quadruples are said to
be {\em join-equivalent} if they give the same image under $\Phi$}
\end{definition}

By straightforward combinatorics, we get that there are
$$|M_Z|!|M_A|!\frac{(|RZ|!)^2}{|\Aut_{\Gamma_A,\Gamma_Z,L}(RZ,RA_\infty)|}$$
graph-join quadruples in $(\Gamma_Z,\Gamma_A,L,J)$'s equivalence
class.


\begin{corollary} \label{totalgluing} {Consider a join equivalence
class of quadruples
$$[\Upsilon]=[(\Gamma_Z,\Gamma_A,L,J)].$$
Let $N=\cm(\ca\sqcup\mz,\Upsilon)$. Let
$$M_\Upsilon=\coprod_{(\Gamma_Z',\Gamma_A',L',J')}
\cm(\ca,\Gamma_A')\times_{D^r}\cm(\cz,\Gamma_Z')$$ where the
disjoint union is over quadruples join-equivalent to $\Upsilon$.
Then $\Phi_{[\Upsilon]}:M\srarr N$ is an \'{e}tale map of degree
$$|M_Z|!|M_A|!(|RZ|!)^2$$}
\end{corollary}

Likewise, we may define graph-join for rubber graphs, $\Gamma_t$,
$\Gamma_b$ (where $t$ and $b$ stand for top and bottom).  Let
$L:R_{b0}\srarr R_{t\infty}$ be a bijective function satisfying
$$\mu_b^0(q)=\mu_t^\infty(L(q)).$$
Let
$$J:M_b\sqcup M_t\srarr \{1,\dots,|M_b|+|M_t|\}$$
be a bijective map.  Then we define the {\em graph join}, a rubber
graph $\Gamma=\Gamma_t*_{L,J} \Gamma_b$ as above, except that
instead of condition (3) above, we have
$$R_0=R_{t0},\ \ a_0=b_{A_t}\circ a_{0t}$$
$$\mu^0=\mu_t^0$$
$$R_\infty=R_{b\infty},\ \ a_\infty=b_{A_b}\circ a_{\infty}b$$
$$\mu^\infty=\mu_b^\infty.$$

Now, let $r=|R_{b0}|=|R_{t\infty}|$.  Exactly as above, we have

\begin{theorem} There is a morphism
$$\Phi:\cm(\ca,\Gamma_{A_t})\times_{D^r} \cm(\ca,\Gamma_{A_b})\srarr \cm(\cz,\Gamma_{A_b}*_{L,J}\Gamma_{A_t}).$$
that is of finite, \'{e}tale of degree
$|\Aut_{\Gamma_{A_b},\Gamma_{A_t},L}(RA_{b0},RA_{t\infty})|$ onto
its image.
\end{theorem}

\begin{corollary} {Consider a moduli stack $N=\cm(\ca,\Gamma_{A_b}
\sqcup_{L,J} \Gamma_{A_t})$.  Let
$$M=\coprod_{(\Gamma_{A_t}',\Gamma_A',L',J')}
\cm(\ca,\Gamma_{A_b}')\times_{D^r}\ca(\ca,\Gamma_{A_t}')$$ where
the disjoint union is of quadruples which give $\Phi$ with image
$N$. The $M\srarr N$ is an \'{e}tale map of degree
$$|M_{A_b}|!|M_{A_t}|!(|R_{A_b0}|!)^2$$}
\end{corollary}

\subsection{Virtual Cycles}

In \cite{Li2}, Li constructs a virtual cycle on $\cm(\cz,\Gamma)$.
This is accomplished by first constructing a perfect obstruction
theory.  This obstruction theory is defined over charts as a
particular two-term complex.  The complex over charts is
constructed by considering charts of
$\cm(Z[n],\Gamma)^{\text{st}}$, a moduli stack of stable maps of
type $\Gamma$ satisfying the pre-deformability and the stability
conditions but not quotiented by the $G[n]=(\cs)^n$-action.  Li
then uses a theorem (\cite{Li2}, Theorem 2.2) to quotient the
complex by the action induced by the Lie algebra of $G[n]$.  This
complex gives a perfect obstruction theory on a chart on
$\cm(\cz,\Gamma)$.  This yields a perfect obstruction theory on
$\cm(\cz,\Gamma)$.

It is completely straightforward to apply this construction to
$\cm(\ca,\Gamma)$.  Instead of using
$\cm(Z[n],\Gamma)^{\text{st}}$, one uses
$\cm(A[n],\Gamma)^{\text{st}}$. One modifies the complex to
enforce particular multiplicities to $D_0$ and $D_\infty$ instead
of $D$. Then one applies Theorem 2.2 of \cite{Li2} to quotient the
complex.  One uses $\Gt{n}$ instead of $G[n]$, but other than
that, the proof is unchanged.  One then obtains a perfect
obstruction theory on $\cm(\ca,\Gamma)$.  The virtual cycle
follows from Li's very general construction.



    \section{Line-Bundle Systems}

We introduce the notion of a line-bundle system.  This is a
sequence of line-bundles $L[n]$ on $\aff^n$, $Z[n]$, or $A[n]$
with certain transition properties.  They will induce line-bundles
on the stack of relative stable and rubber morphisms. In more
esoteric language, line-bundle systems are line-bundles on the
stacks of degenerations $\cz^\rel$, $\ca^\rel$ and their universal
targets. We've chosen to use more primitive notions for the sake
of readability.

\subsection{Definition of line-bundle Systems}

\begin{definition} A {\em base of a line-bundle system} is a
triple $(B[n],H[n],\{i\})$ where

\begin{enumerate}
\item[(1)] $B[n]$ is a sequence of schemes.

\item[(2)] $H[n]$ is a group acting on $B[n]$.

\item[(3)] $\{i\}$, a set of morphisms of swags $i:B[n-1]\srarr
B[n]$ for varying $n$, called standard inclusions.

\end{enumerate}
\end{definition}

The bases that we shall are the following:

\begin{enumerate}
\item[(1)] $\aff=(\aff^n,G[n],\{i\})$ where $\{i\}$ are the
standard inclusions.

\item[(2)] $\cz=(Z[n],G[n],\{i\}).$

\item[(3)] $\bee=(\aff^n,\Gt{n},\{i\}).$

\item[(4)] $\ca=(A[n],\Gt{n},\{i\}).$
\end{enumerate}

\begin{definition} A {\em line-bundle system on
$(B[n],G[n],\{i\})$} is a sequence of $G[n]$-equivariant
line-bundles $L[n]$ on $B[n]$ for each non-negative integer $n$
together with a line-bundle isomorphism $i':L[n-1]\cong i^*L[n]$
for every standard inclusion $i$ such that
\begin{enumerate}
\item[(1)] $L[n]$ is an equivariant line-bundle on $B[n]$ under
group action $H[n]$.

\item[(2)] The group action on $L[n]$ commutes with standard
inclusions. That is, given a standard inclusion
$$i:B[n-1]\srarr B[n]$$
$$i_*: H[n-1]\srarr H[n]$$
then $i^*L[n]$ is isomorphic to $L[n-1]$ as $H[n-1]$-equivariant
bundles.
\end{enumerate}

\end{definition}

\subsection{Induced Line-Bundle Systems}

\begin{proposition}
\begin{enumerate}
\item[(1)] A line-bundle system $L[n]$ on $\aff$ naturally induces
a line-bundle system on $\cz$.

\item[(2)] A line-bundle system $L[n]$ on $\bee$ naturally induces
a line-bundle system on $\ca$. \end{enumerate}
\end{proposition}

\begin{proof} This follows from the fact that the morphisms
$$\begin{array}{ccccc}
p&:&Z[n]&\srarr&\aff^n\\
p&:&A[n]&\srarr&\aff^n.
\end{array}$$
are equivariant and commute with standard inclusions.
\end{proof}

There is a natural notion of a map between line-bundle systems.
Let $B[n]$ be the base of a line-bundle system with group $H[n]$.

\begin{definition} {A map between two line-bundle systems $L[n]$
and $M[n]$ over a base $B[n]$ is a sequence of line-bundle
isomorphisms $f:L[n]\srarr M[n]$ such that
\begin{enumerate}
\item[(1)] $f$ is an equivariant map under $H[n]$

\item[(2)] The map $f$ commutes with standard inclusions.
\end{enumerate}
}
\end{definition}

There is also a natural, obvious notion of tensor product of
line-bundle systems.

We can also define sections of line-bundle systems.

\begin{definition}  A sequence of sections $s[n]:B[n]\srarr L[n]$
is a {\em section of the line-bundle system} if
\begin{enumerate}
\item[(1)] $s[n]$ is equivariant under $H[n]$.

\item[(2)] $s[n]$ commutes with standard inclusions.
\end{enumerate}
\end{definition}

\subsection{Reference Line-Bundle Systems}

We have several a line-bundle system that will be very important
in the sequel, the reference line-bundle systems.

\begin{proposition} {There is a line-bundle system, $V[n]$ on $\cz$
together which a section $s[n]$ that when restricted to $Z[n]$ has
$D[n]$ as its zero divisor.}
\end{proposition}

\begin{proof} $(V[n],s[n])$ is defined by Lemma \ref{pddn}.  It is
straightforward to verify that $V[n]$ is a line-bundle system and
$s[n]$ is a section.
\end{proof}

There is an analogous result for $A[n]$.

\begin{proposition} There are line-bundle systems, $V_0[n],
V_\infty[n]$ together with sections $s_0[n],s_\infty[n]$ on $\ca$
whose corresponding divisor on $A[n]$ are $D_0[n]$, $D_\infty[n]$.
\end{proposition}

\begin{proof}
Let us first construct $V_\infty[n]$, $s_\infty[n]$.

Let $V_\infty[0]=\oh(1)$ with linearization dual to $\oh(-1)$ with
$$\sig{0}\cdot (l,t)\mapsto (\sig{0}l,t).$$
$L_\infty[0]$ has a canonical section $s_\infty[0]$ which is dual
to a section $s_\infty[0]^\vee$ given by
$$s_\infty[0]^\vee:[l:t]\mapsto \left(\frac{l}{t},1\right)$$
$s_\infty[0]$ is an equivariant section with zero scheme
$D_\infty$.  Define $V_\infty[n]$,$s_\infty[n]$ by
\begin{eqnarray*}
V_\infty[n]&=&b^*V_\infty[0]\\
s_\infty[n]&=&b^*s_\infty[0].
\end{eqnarray*}

Likewise $L_0[0]$ to be the line-bundle on $A[0]$ dual to
$\oh(-1)\otimes L^\vee$ with the linearization
$$\sig{0}\cdot ((l,t)\otimes a)=(\sig{0} l, t)\otimes \sig{0}^{-1}
a$$ $s_0[0]$ is defined as the section dual to
$$s_0[0]^\vee:[l:t]\mapsto (l,t)\otimes a_l$$
where $a_l\in L^\vee$ is defined by $a_l(l)=1$.  $s_0[0]$ is an
equivariant section with zero scheme $D_0$.  $V_0[n],s_0[n]$ are
given by Define $V_\infty[n]$,$s_\infty[n]$ by
\begin{eqnarray*}
V_0[n]&=&t^*V_0[0]\\
s_0[n]&=&t^*s_0[0].
\end{eqnarray*}

\end{proof}


    \section{Line-Bundles on the moduli stacks}

In this section, we will define bundles on $\cm(\cz,\Gamma_Z)$ and
$\cm(\ca,\Gamma_A)$. We begin with the bundles on
$\cm(\cz,\Gamma_Z)$. They are
\begin{enumerate}
\item[(1)] $\Dil$, a line-bundle that has a section whose zero
scheme is supported on all split curves. \item[(2)] $\Le{i}$ where
$i$ is an interior marked point (\ref{bodymarkedpoint}), a
line-bundle that has a section whose zero scheme is supported on
split curves where $i$ is mapped to an extended component
(\ref{extendedcomponents}) of the target.
\end{enumerate}

\subsection{Line-bundles on stacks}

Let us recall some definitions from \cite{Vi}.

\begin{definition} Let $F$ be an algebraic stack.  A {\em line-bundle $\cl$} on $F$ can be given by the following data
\begin{enumerate}
\item[(1)] A particular atlas $U$ together with a line-bundle
$\cl_U$ on $U$.

\item[(2)] For the fiber product
$$\xymatrix{
&U\times_F U \ar[dl]_{p_1} \ar[dr]^{p_2}&\\
U&&U}$$
we have an isomorphism $\alpha:p_1^* \cl_U \srarr p_2^*
\cl_U$ that satisfies the cocycle condition.  That is, on
$U\times_F U \times_F U$ with $p_{ij}:U\times_F U \times_F U\srarr
U\times_F U$ ($i,j\in\{1,2,3\}$, $i\neq j$) being projection onto
pairs of factors, we have
$$p_{23}^*\alpha \circ p_{12}^*\alpha = p_{13}^*\alpha:p_1^*\cl_U\srarr
p_3^*\cl_U$$

\end{enumerate}

\end{definition}

We will use atlases which have nice properties over which it will
be simple to define line-bundles.

\begin{definition} \label{pup} A property (P) of a morphism from a
scheme to an algebraic stack $F$ is said to be {\em preserved
under pullback} if given any morphism $S \srarr F$, from a scheme
to $F$ with property $P$ and any morphism of schemes $T \srarr S$
then the composition
$$T\srarr S \srarr F$$
has property (P)
\end{definition}

An example of a property that is preserved under pullback is for a
morphism $S\srarr \cm(\cz,\Gamma)$ to be written as the disjoint
union $S=\coprod S_\alpha$ so that each
$S_\alpha\srarr\cm(\cz,\Gamma)$ is a nice family.  Another example
which we will meet later (Definition \ref{iadmis}) is that of a
family being the disjoint union of $i$-admissible families.

Let us suppose that we have an atlas with a property (P) that is
preserved under pullback.  Now, instead of specifying the
line-bundle on a particular atlas, we can specify it, a fortiori,
over all morphisms $S\srarr F$ with property (P) provided that it
satisfies certain pull-back and transition properties.

\subsection{line-bundles on $\cm(\cz,\Gamma)$ from Line-Bundle Systems}

Now, we restrict to the case where the stack in question is
$\cm(\cz,\Gamma_Z)$. The results are equally true for
$\cm(\ca,\Gamma_A)$ if we replace $Z$'s with $A$'s and $G[n]$'s
with $\Gt{n}$'s. Property (P) will be that a morphism
$S\srarr\cm(\cz,\Gamma)$ can be written as a disjoint union of
nice families
$$\xymatrix{
\ex \ar[r]^{f_\alpha}\ar[d] & Z[n_\alpha]\ar[d]\\
S_\alpha \ar@<1ex>[u]^\gamma \ar[r]_{h_\alpha} & \aff^n}$$

\begin{theorem} {A line-bundle system on $\aff$ induces a line-bundle on $\cm(\cz,\Gamma)$}
\end{theorem}

\begin{proof} As above, we set property (P) to be that for $S\srarr\cm(\cz,\Gamma)$,
$S$ can be written as the disjoint union of $S_\alpha$ where each
$S_\alpha$ is a nice family.

Given a nice family on $\cm(\cz,\Gamma)$,
$$\nicechartonZ{f}{h} $$
define a line-bundle $L$ on $S$ as $h^*L[n]$.  This is
well-behaved under pull-backs.

Now, we have to consider the transitions.  Since the definition
only depends on $h$, not on $f$, we need only consider
isomorphisms $S\times_{\aff^{n_1}} Z[n_1]\srarr S
\times_{\aff^{n_2}} Z[n_2]$ which can be written locally as
effective arrows.  We can factor these as transitions under
standard inclusions and the action of $\rho:S\srarr G[n]$. But
this transition data is exactly specified by the definition of a
line-bundle system.  It is standard to verify that the cocycle
condition is satisfied.
\end{proof}

Given an interior marked point, we can study the line-bundle on
$\cm(\cz,\Gamma)$ induced by a line-bundle system on $\cz$.

\begin{theorem} A line-bundle system on $\cz$ together with the
choice of an interior marked point induces a line-bundle system on
$\cm(\cz,\Gamma)$.
\end{theorem}

\begin{proof} Consider the nice chart
$$\nicechartonZmp{f}{h} $$
where $\gamma$ is the section corresponding to the interior marked
point. We define a bundle $L$ on $S$ by
$$L=(f\circ\gamma)^* L[n].$$
\end{proof}

We have the following simple facts

\begin{lemma} Given line-bundle systems on $\aff$ or $\cz$,
$L[n]$,$M[n]$, $L[n]\otimes M[n]$, then if $L$, $M$, $N$ are the
line-bundles induced by $\cm(\cz,\Gamma)$, respectively then
$$N=L\otimes M$$
\end{lemma}

\begin{lemma} A section of a line-bundle system on $\aff$ or on
$\cz$ gives a section of the induced line-bundle on
$\cm(\cz,\Gamma)$.
\end{lemma}

\subsection{The Dilation Bundle}

In this section, we define the dilation line-bundle, $\Dil$ on
$\cm(\cz,\Gamma)$ induced from a line-bundle system.  The name
"dilation bundle" comes from the dilation of the fibers of
$P=\proj_D(N\oplus 1)$ where $N$ is the normal bundle to $D$ in
$Z$ and is not related to the dilation equation in Gromov-Witten
theory with descendants.

$\Dil$ is induced from a line-bundle system on $\aff$.  We define
the line-bundle system $DIL[n]$ as follows.  Let $DIL[n]$ be the
trivial line-bundle, $1_{\aff^n}$ on $\aff^n$ with $G[n]$ action
given by
$$(\sig{1},\dots,\sig{n})\cdot
s=(\sig{1}\dots\sig{n})^{-1} s$$

\begin{proposition} \label{dilsec} {The line-bundle system $DIL[n]$ has a section.}
\end{proposition}

\begin{proof} Consider the section of $DIL[n]$ over $\aff^n$ defined by
$$s(x_1,\dots,x_n)=x_1 \dots x_n.$$
\end{proof}

$\Dil$ has a simple geometric interpretation which will be proved
later: $c_1(\Dil)$ on $\cm(\cz,\Gamma)$ is (counted with
multiplicity) the locus of split maps.  It is easy to see that is
a reasonable fact. $\Dil$ has a section whose zero scheme are
split maps: this section is induced from a section of the
line-bundle system $DIL[n]$ on $\aff^n$; the zero-scheme of this
section is the the union of hyper-planes in $\aff^n$ consisting of
points with at least one coordinate equal to $0$; the fiber of
$Z[n]\srarr \aff^n$ over any closed point in this scheme consists
of a non-smooth target.

\subsection{Definition of $\Le{i}$}

Consider $\cm(\cz,\Gamma)$ with at one distinguished interior
marked point, which we will call $i$.  We define a line-bundle
$\Le{i}$ which will have a section whose zero stack is supported
on split maps where the marked point $i$ is mapped to an extended
component of the target. This line-bundle is induced from a
line-bundle system on $\cz$.

We define $LE[n]$ on $Z[n]$.  Recall that $Z[0]=Z$. $D\subset Z$
is a Cartier divisor and is the zero-section of a section $s$ of a
line-bundle $V$ on $Z$.  Define $LE[0]=V$, $LE[n]=c^*LE[0]$ where
$c$ is the collapsing map
$$c:Z[n]\srarr Z.$$

\begin{lemma} The line-bundle system $LE[n]$ has a section.
\end{lemma}

\begin{proof} Define the section by
$$se[n]=c^*s.$$
\end{proof}

Note that if we take a closed point $p\in\aff^n$ then $s^*c$ is
zero on the extended components of $Z[n]\times_{\aff^n} p$.

\begin{theorem} $\Le{i}$ is canonically isomorphic to $\ev_i^*L$.
\end{theorem}

\begin{proof} On $S\srarr\cm(\cz,\Gamma)$, a family of $\cm(\cz,\Gamma)$, we have the
diagram
$$\xymatrix{
\ex \ar[d] \ar[r]^f & Z[n]\ar[d]\ar[r]^c&Z\\
S \ar@<1ex>[u]^{\gamma}\ar[r] & \aff^n&}$$

But, $\ev_i^*L=(c\circ f \circ \gamma)^*L$ which is canonically
the line-bundle induced on $S$ by the line-bundle system $LE[n]$.
\end{proof}

\subsection{Line-Bundles on $\cm(\ca,\Gamma)$ from Line-Bundle Systems}

In this subsection, we define the following bundles on
$\cm(\ca,\Gamma)$ which are induced from line-bundle systems:
\begin{enumerate}
\item[(1)] $\Top$, the target cotangent line-bundle at $X_0$.
$c_1(\Top)=\Psi_0$ is the target $\Psi$ class of \cite{FP}.

\item[(2)] $\Bot$, the target cotangent line-bundle at $X_\infty$
which is $\Top$'s upside-down analog. $c_1(\Bot)=\Psi_\infty$.

\item[(3)] $\Split$, the Split bundle which has a section whose
zero stack is supported on all split maps (Definition
\ref{splitmaps}).

\item[(4)] $\Lnt{i}$, the not-top bundle with respect to a
distinguished interior marked point $i$.  This bundle has a
section whose zero stack is supported on split maps where the
$i$th marked point is not on a top component.

\item[(5)] $\Lnb{i}$, the not-bottom bundle with respect to a
interior marked point $i$.  This bundle has a section whose zero
stack is supported on split maps where the $i$th marked point is
not on a bottom component.
\end{enumerate}

We prove the following relations among bundles which will be used
in the sequel.
\begin{eqnarray*}
\Top\otimes\Bot&=&\Split\\
\Top\otimes\ev_i^*L^\vee&=&\Lnt{i}\\
\Bot\otimes\ev_i^*L&=&\Lnt{i}
\end{eqnarray*}

\subsection{Line-bundles on $\cm(\ca,\Gamma)$ from line-bundle systems}

Line-bundles can be defined on $\cm(\ca,\Gamma)$ by defining them
over nice families, or more generally on charts with a property
(P) preserved under pullback (\ref{pup}) so that there is an atlas
with property (P).

\subsection{Definition of $\Top$ and $\Bot$}

$\Top$ is induced from a line-bundle system on $\bee$.

We define the line-bundle system $T[n]$ as follows.  Let $T[n]$ be
the trivial line-bundle, $1_{\aff^n}$ on $\aff^n$ with $\Gt{n}$
action given by
$$(\sig{0},\sig{1},\dots,\sig{n})\cdot
s=(\sig{0}\sig{1}\dots\sig{n})^{-1} s$$ Therefore the maps in the
following diagram

$$\xymatrix{
T[n]\ar[r]^{\sigma \cdot}\ar[d]& T[n]\ar[d]\\
\aff^n \ar[r]&\aff^n}$$
are given as
$$\begin{array}{rcccl}
& s & \mapsto & (\sig{0}\sig{1}\dots\sig{n})^{-1} s &\\
&\downarrow & & \downarrow & \\
& (x_1,\dots,x_n) & \mapsto & (\sig{1}^{-1}x_1,\dots,\sig{n}^{-1} x_n)& \\
\end{array}$$

The transition map under an effective inclusions $i$ is trivial.
We call the induced line-bundle $\Top$.

The definition of $\Bot$ is similar to that of $\Top$. It is
induced by a line-bundle system $B[n]$. $B[n]$ is the trivial
line-bundle $1_{\aff^n}$ on $\aff^n$ under the $\Gt{n}$-action
$$(\sig{0},\sig{1},\dots,\sig{n})\cdot s=\sig{0} s$$
and trivial transition map under effective inclusion.

$\Top$ and $\Bot$ can be given an interpretation in the stack of
rational sausages, the substack of $\acm_{0,2}$ consisting of
pre-stable curves so that the two marked points lie on different
sides of every node. $\Top$ and $\Bot$ are equal to the pullbacks
of the cotangent line classes at the two marked points. See
\cite{GV} for an elaboration.

\subsection{The split bundle}

We define $\Split$ to be the bundle on $\cm(\ca,\Gamma)$ induced
from the line-bundle system on $\bee$ given by the trivial bundle
$S[n]=1_{\aff^n}$ on $\aff^n$ with the trivial transition map
under effective inclusion and the following $\Gt{n}$-action:
$$(\sig{0},\sig{1},\dots,\sig{n})\cdot
s=(\sig{1}\dots\sig{n})^{-1} s.$$

\begin{proposition} {The line-bundle system $S[n]$ has a section.}
\end{proposition}

\begin{proof} The proof is identical to the one for $DIL[n]$,
(Proposition \ref{dilsec}). \end{proof}

\begin{proposition} $\Top\otimes\Bot=\Split$
\end{proposition}

\begin{proof} Because $T[n]\otimes B[n]=S[n]$ as line-bundle
systems, the induced bundles are equal. \end{proof}

$\Split$ has a simple geometric interpretation analogous to that
of $\Dil$.  On $\aff^n$ under $\Gt{n}$, $\Split$ has a section
given by $x_1x_2\dots x_n$.  The zero-scheme of this section is
the the union of hyper-planes in $\aff^n$ consisting of points
with at least one coordinate equal to $0$.  We know that the fiber
of $A[n]\srarr \aff^n$ over any closed point in this scheme
consists of a chain of more than one $P$'s.  These sections glue
together to form a section of $\Split$ on $\cm(\ca,\Gamma)$. It
follows that $c_1(\Split)$ on $\cm(\ca,\Gamma)$ is (counted with
multiplicity) the locus of split maps.

If we consider the stack of rational sausages where $\Top$ and
$\Bot$ are the restriction of $\psi$ classes on $\acm_{0,2}$, this
relation is the pullback of the genus 0 recursion relation of Lee
and Pandharipande \cite{LP}.

\subsection{Definition of $\Lnt{i}$}

Consider $\cm(\ca,\Gamma)$ a moduli stack of rubber maps with a
distinguished interior marked point, $i$.  We will define two
bundles $\Lnt{i}$ and $\Lnb{i}$ on $\cm(\ca,\Gamma)$ with very
clear geometric meanings. $c_1(\Lnt{i})$ will be supported on the
substack of relative stable maps consisting of split maps where
the $i$th interior marked point does not lie on the top component
of the target while $c_1(\Lnb{i})$ will be the substack where the
$i$th interior marked point does not lie on the bottom component.
Many geometric results follow from equations relating $\Lnt{i},
\Lnb{i}$ to other line-bundles.

$\Lnt{i}$ and $\Lnb{i}$ are induced from line-bundle systems on
$\ca$, $NT[n]$ and $NB[n]$, respectively.

Let $NT[0]$ be the equivariant bundle on
$$A[0]=\proj(L\oplus 1)$$
which is dual to $\oh(-1)$ with the $\Gt{0}$-action
$$\sig{0}\cdot(l,t)=(\sig{0}l,t).$$
Let $NT[n]=t^*NT[0]$ where $t:A[n]\srarr A[0]$ is the map given in
Definition \ref{topmap}

Note also that $NT[0]$ has a section $s$ whose zero-scheme is
$D_\infty$ so $NT[n]$ has a section $t^*s$. If we look $_l P$, the
fiber over a point $x\in\aff^n$, we see that this section is zero
on $P_0,\dots,P_{l-1}$ and is non-zero on $P_l$ away from
$D_\infty$.

For $i$, an interior marked point for $\cm(\ca,\Gamma)$, $\Lnt{i}$
is the line-bundle on $\cm(\ca,\Gamma)$ induced by $NT[n]$ and
$i$.

$\Lnb{i}$ is defined similarly.  Let $NB[0]$ be the equivariant
bundle on $A[0]=\proj_X(L\oplus 1)$ dual to $\oh(-1)\otimes
L^\vee$ with the $\Gt{0}$-action-
$$\sig{0}\cdot((l,t)\otimes a)=(\sig{0}l,t)\otimes \sig{0}^{-1}a.$$
$NB[n]$ is given by $b^*NB[0]$.

For $i$, an interior marked point for $\cm(\ca,\Gamma)$, $\Lnb{i}$
is the line-bundle on $\cm(\ca,\Gamma)$ induced by $NB[n]$ and
$i$.

\subsection{The Reference Line-Bundles}
Earlier, we defined line-bundle systems $V_0[n]$, $V_\infty[n]$ on
$\ca$.

\begin{proposition} \label{d0trivial} {The line-bundles induced by
$V_0[n]$, $V_\infty[n]$ on $\cm(\ca,\Gamma)$ are trivial}
\end{proposition}

\begin{proof} We produce a canonical nowhere-zero section of the
induced line-bundle by $V_0[n]$.  The case of $V_\infty[n]$ is
analogous.  Recall that the line-bundles are defined on each nice
chart
$${\xymatrix{
\ex \ar[d] \ar[r]^f & A[n]\ar[d]\ar[r]^t&A[0]\\
S \ar@<1ex>[u]^{\gamma}\ar[r] & \aff^n&}}$$

$V_0$ is defined on $S$ to be
$$V_0=(f\circ \gamma)^*V_0[n].$$
It has a section $(f\circ \gamma)^*s_0[n]$.  Since $s_0[n]$ is
non-zero away from $D_0[n]$ and interior marked points are not
mapped to $D_0[n]$, the pulled-back section is non-zero.  It is
easy to verify that this section transformers properly.
\end{proof}

\subsection{The Evaluation Map and Line-bundle System}

\begin{lemma}
Let $i$ be an interior marked point on $\cm(\ca,\Gamma)$, and let
$M$ be a line-bundle on $X$.  Then there is a line-bundle system
$M[n]$ on $\ca$ so that the line-bundle it induces on
$\cm(\ca,\Gamma)$ at $i$ is $\ev_i^*M$
\end{lemma}

\begin{proof} Given a nice chart, we have the diagram

$${\xymatrix{
\ex \ar[d] \ar[r]^f & A[n]\ar[d]\ar[r]^{\pi}&X\\
S \ar@<1ex>[u]^{\gamma}\ar[r] & \aff^n&}}$$

Let $M[n]=\pi^*M$.  $M[n]$ is easily seen to be a line-bundle
system, and since $\pi\circ f \circ \gamma=\ev_i$, induced
line-bundle is $\ev_i^*L$.
\end{proof}

\subsection{The Marked Point Interpretation of the Top Bundle}

By relating the various bundles described in this section, we can
find significant geometric facts.  There is one geometric
interpretation of the top bundle that will be supremely useful in
the sequel.  This interpretation was first advanced by Andreas
Gathmann in conversation with the author.

Let $\cm(\ca,\Gamma)$ be the moduli stack of rubber maps with a
distinguished interior marked point, $i$.  Consider the evaluation
map
$$\ev_i:\cm(\ca,\Gamma)\srarr X$$ at $i$.

\begin{theorem} \label{mpiotb} $\Top\otimes \ev_i^* L^\vee =
\Lnt{i}$
\end{theorem}

\begin{proof} It suffices to prove that $\Top^\vee \otimes \ev_i^* L
\otimes \Lnt{i}$ is the trivial bundle.  To do this, we will show
that the line-bundle system on $\ca$ that induces $\Top^\vee
\otimes \ev_i^* L \otimes \Lnt{i}$ is line-bundle
system-isomorphic to $V_0[n]$ which induces the trivial bundle on
$\cm(\ca,\Gamma)$.

Consider the following diagram of swags where $X$ is equipped with
the trivial group action
$$\xymatrix{
A[n]\ar[dr]_\pi\ar[rr]&&A[0]\ar[dl]^\pi\\
&X&}$$
Consider also the map of swags $p:A[n]\srarr \aff^n$.

Let us review the bundles involved  They are pullbacks of bundles
from $A[0]$ by $t^*$.
\begin{enumerate}
\item[(1)] $\Top^\vee$ is induced by the bundle
$\theta[n]=p^*T[n]^\vee$ on $A[n]$, where $T[n]^\vee$ is the
trivial bundle on $\aff^n$ with the $\Gt{n}$-linearization
$$(\sig{0},\sig{1},\dots,\sig{n})\cdot s =
(\sig{0}\sig{1}\dots\sig{n}) s.$$ Note that
$\theta[n]=t^*\theta[0]$ where $\theta[0]$ is the trivial bundle
on $A[0]$ under the action
$$\sig{0}\cdot s= \sig{0} s.$$

\item[(2)] $\ev_i^* L$ is induced by the line-bundle system
$\pi^*L$ on $A[n]$.  But $\pi^* L=t^* \pi^* L$.

\item[(3)] $\Lnt{i}$ is induced by $NT[n]=t^*NT[0]$.

\item[(4)] A trivial bundle is induced by $L_0[n]=t^*L_0[0]$.
\end{enumerate}

Therefore, to prove $\theta[n]\otimes e^*L \otimes NT[n]=L_0[n]$,
as a line-bundle system on $A[n]$, it suffices to verify
$$\theta[0]\otimes e_0^* L \otimes NT[0]=L_0[0]$$
on $A[0]=P$ which is straightforward.
\end{proof}

Likewise,
\begin{theorem} \label{mpiobb} {$\Bot\otimes \ev_i^* L = \Lnb{i}$}
\end{theorem}

\begin{proof} The proof is exactly analogous to the above. \end{proof}

\subsection{Boundary $\Psi$ class interpretation of $\Top$}

There is an interpretation of $\Top$ in terms of $\psi$ classes at
boundary marked points.  Consider $\cm(\ca,\Gamma)$.  Let $j$ be a
boundary marked point evaluating to $X_0$ with multiplicity $m$.
Let $L_j$ be the tangent space at the $j$th marked point. Then, we
will show $\Top=(L_j^\vee)^{m} \otimes \ev_i^*L$.  This relation
tells us that nothing new can be found in rubber theory with
$\psi$ classes at the boundary marked point. Also, since this
formula is independent of $i$, we can relate $\psi$ classes at
different boundary marked points.

\begin{theorem} {$\Top=(L_j^\vee)^{\otimes m}\otimes \ev_j^* L$}
\end{theorem}

\begin{proof} We prove this by exhibiting a regular, nowhere vanishing
section of this bundle.

We need to consider the line-bundle on $\cm(\ca,\Gamma)$ induced
from a line-bundle system on $\ca$ by a boundary marked point.
This is identical to a line-bundle induced by an interior marked
point except for the following modification.  Given a line-bundle
system $L[n]$, consider the nice chart on $\cm(\ca,\Gamma)$
$$\xymatrix{
\ex \ar[d] \ar[r]^f & A[n]\ar[d]\\
S \ar@<1ex>[u]^{\delta}\ar[r] & \aff^n}$$ where $\delta$ is the
section corresponding to the boundary marked point evaluating to
$D_0$.  Then, set $L=(f\circ\delta)^*L[n]$. This is seen to give a
line-bundle on $\cm(\ca,\Gamma)$.  Let $V$ be the line-bundle
induced by $V_0[n]$ at $j$.

\begin{lemma} \label{Latbdrymp} {$V=\Top^\vee\otimes\ev_j^* L$}
\end{lemma}

\begin{proof} We show that $V_0[n]^\vee\otimes p^*T[n]^\vee \otimes \pi^*
L$ has an equivariant section that does not vanish near the image
of $(f\circ \gamma)$.

Let us recall the morphism of swags
$$\begin{array}{ccccc}
t&:&A[n]&\srarr&A[0]\\
p&:&A[n]&\srarr&\aff^n\\
\pi&:&A[n]&\srarr&X
\end{array}$$

Note that $\Top^\vee$ is induced at the boundary marked point $i$
from the line-bundle system $p^*T[n]^\vee$. Note also that
$V_0[n]^\vee=t^*V_0[0]^\vee$, $\pi^* L=t^* \pi^* L$, and
$p^*T[n]^\vee=t^*p^*T[0]^\vee$. Therefore, it suffices to show
$L_0[0]^\vee\otimes p^*T[0]^\vee \otimes \pi^* L$ has an
equivariant section that does not vanish near $D_0$.

Now, $L_0[0]^\vee$ is the bundle on $A[0]$ equal to
$\oh(-1)\otimes \pi_0^* L^\vee$ with the linearization
$$\sig{0}\cdot ((l,t)\otimes a)=(\sig{0} l, t)\otimes \sig{0}^{-1} a.$$

It follows that $L_0[0]^\vee\otimes p_0^*T[0]^\vee \otimes \pi_0^*
L$ is $\oh(-1)$ with the linearization
$$\sig{0}\cdot ((l,t)) = (\sig{0} l, t)$$
But $\oh(-1)$ has an invariant rational section that is
well-defined and nonzero near $D_0=\{[0:t]\}$.
\end{proof}

\begin{lemma} There is a nowhere vanishing regular section of $(L_j)^{\otimes m}\otimes V$.
\end{lemma}

\begin{proof}
We observe that on a nice chart,
$$\delta^*(\Omega_{\ex/S}^{\otimes m} \otimes f^*V_0[n])$$
has a $\Gt{n}^S$-equivariant, nowhere vanishing, regular section.
Since $\ex/S$ is smooth near $\gamma(S)$, we have
$$\delta^*(\Omega_{\ex/S})^{\otimes m} \otimes f^*V_0[n]=
\Shom(N^{\otimes m}_{\delta(S)/\ex},N_{D_0[n]/A[n]}).$$ But this
bundle has a section given by the projection of the $m$th formal
derivative of $f$ to $N_{D_0[n]/A[n]}$.
\end{proof}

\end{proof}




    \section{Degeneration Formulae}

In the previous section, we have proved equations relating the
line-bundle
$$\Dil,\Split,\Lnt{i},\Lnb{i}$$
to other line bundles on $\cm(\cz,\Gamma_Z)$ and
$\cm(\ca,\Gamma_A)$. The first Chern classes of these line-bundles
represent specific geometric situations involving split curves.
For example, $c_1(\Split)$ is a substack of $\cm(\ca,\Gamma)$,
that is, in a virtual sense, all split curves.  $c_1(\Lnt{i})$
virtually consists of all split curves in which the $i$th marked
point is not on the topmost component.  Of course, we are counting
certain curves with multiplicity and there are virtual issues that
make geometric interpretations inaccurate. Fortunately, by closely
adapting the arguments in \cite{Li2}, we can make arguments
compatible with virtual cycle constructions. This allows us to
write the cap product of a first Chern class of one of our bundles
with the virtual cycle in terms of the virtual cycles of smaller
moduli spaces.

We will express the first Chern class of various line-bundles
geometrically by adapting Li's argument \cite{Li2}. The argument
is in several stages and we state it only in the case
$\cm(\cz,\Gamma)$ noting that the case for $\cm(\ca,\Gamma)$ is
exactly analogous:
\begin{enumerate}
\item[(1)] For $\Gamma$, consider quadruples
$\Upsilon=(\Gamma_Z,\Gamma_A,L,J)$ so that the graph join,
$\Gamma_Z*_{L,J}\Gamma_A$ is isomorphic to $\Gamma$ We can define
a line bundle $L_\Upsilon$ on $\cm(\cz,\Gamma_Z)$.

\item[(2)] We show that
$$c_1(L_\Upsilon)\cap \vir{\cm(\cz,\Gamma)}=m(\Upsilon)\vir{\cm(\ca\sqcup\cz,\Upsilon)}$$
where $m(\Upsilon)$ is defined in Definition \ref{multijoin} and
$\vir{\cm(\ca\sqcup\cz,\Upsilon)}$ is an appropriately defined
virtual cycle.

\item[(3)] Given the joining morphism
$$\Phi:\cm(\ca,\Gamma_A)\times_{D^r}
\cm(\cz,\Gamma_Z)\srarr\cm(\cz,\Gamma_A\sqcup_{L,J}\Gamma_Z)$$ and
the diagram
$$\xymatrix{\cm(\ca,\Gamma_A)\times_{D^r}
\cm(\cz,\Gamma_Z)\ar[r] \ar[d]
&\cm(\ca,\Gamma_A)\times \cm(\cz,\Gamma_Z)\ar[d]\\
D^r\ar[r]_\Delta & D^r\times D^r}$$ where $\Delta$ is the diagonal
map. We have
$$\Phi_*\Delta^!(\vir{\cm(\ca,\Gamma_A)}\times
\vir{\cm(\cz,\Gamma_Z)})=\vir{\cm(\ca\sqcup\cz,\Upsilon)}.$$

\item[(4)] Given a line-bundle $L=\Dil$ or $L=\Le{i}$ (or in the
case of $\cm(\ca,\Gamma)$,
$\Split$,$\Lnt{i}$,$\Lnb{i}$), we
exhibit a set of join-equivalence classes $\Omega$ so that
$$L=\otimes_{[\Upsilon]\in\Omega} L_\Upsilon$$

\item[(5)] Consequently
$$c_1(L) \cap \vir{\cm(\cz,\Gamma)}=\sum_{\Upsilon\in\Omega}
m(\Upsilon)\Phi_*\Delta^!(\vir{\cm(\ca,\Gamma_A)}\times\vir{\cm(\cz,\Gamma_Z)})$$
\end{enumerate}

To modify this argument to work for $\cm(\cz,\Gamma)$, replace all
pairs $(\Gamma_A,\Gamma_Z)$ with $(\Gamma_t,\Gamma_b)$ and replace
$\cz$ with $\ca$. (4) is the only item significantly different
from \cite{Li2} to warrant much explanation.

\subsection{Local Interpretation of $\Lnb{i}$}

We study a section of $\Lnb{i}$ in the interest of proving (4)
above.  Recall that $\Lnb{i}$ is induced from a line bundle system
on $\ca$ called $NB[n]$. For the bottom map $b:A[n]\srarr A[0]$,
$NB[n]=b^*NB[0]$ where $NB[0]=V_0[0]$. Therefore, $NB[0]$ is given
on $A[0]=\proj_X(L\oplus 1)$ by a bundle dual to $\oh(-1)\otimes
L^\vee$ under the linearization
$$\sig{0}\cdot((l,t)\otimes a)=((\sig{0}l,t)\otimes
\sig{0}^{-1}a).$$ Therefore, $NB[0]$ has an equivariant section
$s[0]$ given by
$$s[0]^\vee:[l:t]\srarr (l,t)\otimes a_l$$
where $a_l\in L^\vee$ satisfies $a_l(l)=1$.

Now, $s[n]=b^*s[0]$ forms a line bundle system section of $NB[n]$
and therefore induces a section of $\Lnb{i}$ on $\cm(\ca,\Gamma)$.

We need to describe how $s$ looks in charts.

Recall that $K_k$ is the hyperplane on $\aff^n$ cut out by $t_k=0$
and that
$$A[n]\times_{\aff^n} K_k = (A[k-1]\times \aff^{n-k}) \sqcup_D
(\aff^{k-1}\times A[n-k])$$

\begin{definition} \label{iadmis} A chart with a marked point,
$i$,
$$\nicechartonAmp{f}{h} $$
is said to be {\em $i$-admissible for $l\in\{0,1,\dots,n\}$} if
\begin{enumerate}
\item[(1)] if $l\geq 1$ then $S\times K_l$ is nonempty

\item[(2)] {For $k\in\{1,2,\dots,n\}$ so that $S\times_{\aff^n}
K_k$ is non-empty, the image of every closed point under
$$f\circ \gamma:S\times_{\aff^n} K_k\srarr A[n]\times_{\aff^n} K_k=
(A[k-1]\times \aff^{n-k}) \sqcup_D (\aff^{k-1}\times A[n-k])$$
lies in
\begin{enumerate}
\item[(a)] $\aff^{k-1}\times A[n-k]$  for $k\leq l$

\item[(b)] $A[k-1]\times \aff^{n-k}$ for $k\geq l+1$.
\end{enumerate}
}
\end{enumerate}
\end{definition}

The property of $i$-admissibility has a simple geometric
explanation. Suppose a chart is $i$-admissible for
$l\in\{0,1,\dots,n\}$.  Let $s\in S$ be a closed point, $t=h(s)$.
Let $\{a(1),a(2),\dots,a(m)\}$ be a subset of $\{1,\dots,n\}$
corresponding to the coordinates of $t$ which are zero.  We have a
morphism
$$f:C=\ex\times_S s\srarr A[n]\times_{\aff^n} t$$
Then $A[n]\times_{\aff^n} t=\ _m
P=P_0\sqcup_{D_1}\dots\sqcup_{D_m} P_m$. The $i$th marked point is
mapped into the component of $P_i\setminus D_{a^{-1}(l)}$ that
contains $X_0$.

It is easy to see that the property of a morphism
$S\srarr\cm(\ca,\Gamma)$ being $i$-admissible is preserved under
pullbacks in the sense of \ref{pup}

\begin{lemma} {Give a chart
$$\xymatrix{
\ex \ar[d] \ar[r]^f & A[n]\ar[d]\\
S \ar@<1ex>[u]^{\gamma}\ar[r]_h & \aff^n}$$
and a closed point
$p\in S$, then there is a Zariski neighborhood $U_p\subseteq S$ of
$p$ so that
$$\xymatrix{
\ex\times_S U_p \ar[d] \ar[r]^f & A[n]\ar[d]\\
U_p \ar@<1ex>[u]^{\gamma}\ar[r]_h & \aff^n}$$ is $i$-admissible.}
\end{lemma}

\begin{proof} The family over $p$,
$$\xymatrix{
\ex\times_S p \ar[d] \ar[r]^f & A[n]\ar[d]\\
p \ar@<1ex>[u]^{\gamma}\ar[r]_h & \aff^n}$$
is $i$-admissible for
some $l\in\{0,1,\dots,n\}$.

Consequently, we can find a neighborhood $U_k$ of $p$ so that
image of closed points under
$$f\circ \gamma:U_k\times_{\aff^n} K_k\srarr A[n]\times_{\aff^n}=
(A[k-1]\times \aff^{n-k}) \sqcup_D (\aff^{k-1}\times A[n-k])$$
lies lies in
\begin{enumerate}
\item[(1)] $\aff^{k-1}\times A[n-k]$ for $k\leq l$

\item[(2)] $A[k-1]\times \aff^{n-k}$ for $k\geq l+1$.
\end{enumerate}

Let $U_p=U_1\cap U_2\cap \dots \cap U_n$. \end{proof}

\begin{corollary} $\cm(\ca,\Gamma)$ has \'{e}tale atlas $U\srarr\cm(\ca,\Gamma)$ that has
the property that $U$ can be written as a disjoint union,
$U=\coprod U_\alpha$ where each $U_\alpha$ is $i$-admissible.
\end{corollary}

\begin{proof} We can refine an atlas of $\cm(\ca,\Gamma)$ that is the disjoint union of nice charts.
\end{proof}
Note that if
$$\nicechartonAmp{f}{h} $$
is $i$-admissible for an $l\in\{1,\dots,n\}$ and if $j:[n]\srarr
[N]$ is an order-preserving inclusion, then if $j_*:A[n]\srarr
A[N]$, $J:\aff^n\srarr \aff^N$ are the induced inclusions, then
$$\nicechartonAmp{j_*\circ f}{J\circ h} $$
is $i$-admissible for $j(l)$ if $l\geq 1$ and $i$-admissible for
$0$ otherwise.

Note also that the property of $i$-admissibility for
$l\in\{0,1,\dots,n\}$ is invariant under the action of
$\rho:S\srarr \Gt{n}$.

\begin{lemma} \label{localmodel} {Let $S$ be an $i$-admissible chart for
$l\in\{0,1,\dots,n\}$. Then $\Lnb{i}|_S\cong 1_S$, the
topologically trivial bundle under the $\Gt{n}$-action
$$(\sig{0},\sig{1},\dots,\sig{n})[s]=[(\sig{1}\dots\sig{l})^{-1}s]$$
Furthermore, there is a constant $c$ so that the section is given
on $S$ as $h^*(cx_1\dots x_l)$}
\end{lemma}

\begin{proof} We need to define $\Gt{n}$-equivariant line bundles $C_i[n]$
on $\aff^n$.  Let $C_i[n]$ as $1_{\aff^n}$, the topologically
trivial bundle on $\aff^n$ with linearization
$$(\sig{0},\sig{1},\dots,\sig{n})\cdot s = \sig{i}^{-1} s$$
$C_i$ has a $\Gt{n}$-equivariant section given by $x_i$.

We will show that given an $i$-admissible atlas on
$\cm(\ca,\Gamma)$, $\{\ex_\alpha, S_\alpha,
f_\alpha,h_\alpha,\gamma_\alpha\}$, so that
$$\xymatrix{
\ex_\alpha \ar[d] \ar[r]^f & A[n_\alpha]\ar[d]\\
S_\alpha \ar@<1ex>[u]^{{\gamma}_\alpha}\ar[r]_h &
\aff^{n_\alpha}}$$
is $i$-admissible for
$l_\alpha\in\{0,1,\dots,n_\alpha\}$ then
$$(f_\alpha\circ \gamma_\alpha)^*NB[n]\otimes
h_\alpha^*(C_1[n]^\vee\otimes\dots\otimes C_{l_\alpha}[n]^\vee)$$
is canonically isomorphic to the trivial bundle on $S_\alpha$.
 Therefore, we conclude that
$$h_\alpha^*(C_1[n]\otimes\dots\otimes
C_{l_\alpha}[n])$$ induces a line-bundle on $\cm(\ca,\Gamma)$ that
is isomorphic to $\Lnb{i}$.  The following lemma gives a
description of this bundle.

\begin{lemma} \label{lbiad} {There is line bundle on $\cm(\ca,\Gamma)$ defined
over an $i$-admissible chart $S_\alpha\srarr\cm(\ca,\Gamma)$ by
$$h_\alpha^*(C_1[n]\otimes\dots\otimes
C_{l_\alpha}[n])$$ which has a section
$$h_\alpha^*(x_1 x_2\dots x_{l_\alpha})$$}
\end{lemma}

\begin{proof} We just have to check that the line-bundles are well-behaved
under the group action and standard inclusions.  This is standard
in light of the observation preceding Lemma
\ref{localmodel}.\end{proof}

Recall that $\Gt{n}^S$ is the group of morphisms $\rho:S\srarr
\Gt{n}$ under point-wise multiplication.  To obtain the triviality
of $(f_\alpha\circ \gamma_\alpha)^*NB[n]\otimes
h_\alpha^*(C_1[n]^\vee\otimes\dots\otimes C_{l_\alpha}[n]^\vee)$,
we have the following.

\begin{claim} {For any family
$$\nicechartonAmp{f}{\ } $$
that is $i$-admissible for $l\in\{0,1,\dots,n\}$, the bundle
$$(f\circ\gamma)^*(NB[n]\otimes
C_1[n]^\vee\otimes\dots\otimes C_l[n]^\vee)$$
has a nowhere zero
section $u[n]$ that obeys the following properties:
\begin{enumerate}
\item[(1)] $u[n]$ is $\Gt{n}^S$-equivariant

\item[(2)] Given an effective inclusion induced by
$j:[n]\hookrightarrow [N]$, $\rho:S\srarr \Gt{N}$ then there is an
isomorphism
\begin{eqnarray*}
(j_*\circ f\circ\gamma)^*(NB[N])&\otimes &(J\circ
h)^*(C_1[N]^\vee\otimes\dots\otimes
C_{j(l)}[N]^\vee)\\
&\cong&(f\circ \gamma)^*(NB[n]\otimes
C_1[n]^\vee\otimes\dots\otimes C_l[n]^\vee) \end{eqnarray*} that
is $\Gt{n}^S=(j_*\Gt{n})^S$-equivariant and that takes $u[N]$ to
$u[n]$

\item[(3)] Given a morphism $i:T\srarr S$ corresponding to pulling
back charts, $i^*u[n]=u[n]$.

\end{enumerate}}
\end{claim}

We first relate $NB[n]$ to $V_0[n]=t^*V_0[0]$ which has a
canonical equivariant section $s[n]=t^*s[0]$. We note that the
bottom map can be factored as a morphism of swags into blow-downs
and projections as follows:
$$b:A[n]\srarr A[n-1]\times\aff^1\srarr
A[n-1]\srarr\dots\srarr A[1]\srarr A[0]\times\aff^1\srarr A[0]$$
so if
$$\xymatrix{
\pi_n:A[n]\ar[r]^>>>>>\beta & A[n-1]\times\aff^1 \ar[r]^>>>>>q
&A[n-1]}$$ we have
$$NB[n]=\pi_n^* NB[n-1]=\beta^*q^* NB[n-1]$$
On the other hand, $V_0[n]$ is the proper transform of
$q^*V_0[n-1]$.  Consequently, if $E_n=\oh(1)$ is the line bundle
whose zero scheme is exceptional divisor of $\beta$ then
$$V_0[n]=\pi_n^* V_0[n-1]\otimes E_n^\vee.$$ Let
$\pi_{n,i}:A[n]\srarr A[i]$ be given by the composition
$$\pi_{n,i}=\pi_{i+1} \circ \dots \circ \pi_n :A[n]\srarr
A[n-1]\dots A[i+1] \srarr A[i]$$ with the group homomorphism given
by
$$\pi_{n,i*}:(\sig{0},\sig{1},\dots,\sig{n})\srarr
(\sig{0},\sig{1},\dots,\sig{i}).$$ Then, it is easy to see by
induction that we have the following $\Gt{n}$-equivariant
isomorphism of line-bundles,

$$NB[n]=V_0[n]\otimes \pi_{n,1}^* E_1\otimes \pi_{n,2}^* E_2 \otimes
\pi_{n,n-1}^* E_{n-1} \otimes E_n$$

Since $V_0[n]$ induces the trivial bundle on $\cm(\ca,\Gamma)$,
$NB[n]$ induces the same bundle as
$$NB[n]\otimes L_0[n]^\vee=\pi_{n,1}^* E_1\otimes \pi_{n,2}^* E_2 \otimes
\pi_{n,n-1}^* E_{n-1} \otimes E_n.$$

The theorem then follows as a consequence of the following lemma.

\begin{claim}
Consider an $i$-admissible atlas as above.  On each $S_\alpha$,
consider the bundle
\begin{eqnarray*}
(f\circ\gamma_\alpha)^*(\pi_{n_\alpha,1}^* E_1\otimes
\pi_{n_\alpha,2}^* E_2 \otimes \dots\otimes
\pi_{n_\alpha,n_\alpha-1}^* E_{n_\alpha-1} \otimes E_{n_\alpha})\\
\ \ \otimes h_\alpha^*(C_1[n]^\vee \otimes C_2[n]^\vee \otimes
\dots \otimes C_{l_\alpha}[n]^\vee)
\end{eqnarray*}
on $S_\alpha$. These bundles descend to a section to
$\cm(\ca,\Gamma)$. Moreover, they possess a nowhere zero section
that descends to a section on $\cm(\ca,\Gamma)$.
\end{claim}

This claim follows from examining blow-ups in local coordinates.

\end{proof}

There is the analog for $\Lnt{i}$ which has a similarly defined
canonical section.

\begin{lemma} {Let $S$ be an $i$-admissible chart for
$l\in\{0,1,\dots,n\}$. Then $$\Lnt{i}|_S\cong 1_S,$$ the
topologically trivial bundle under the $\Gt{n}$-action
$$(\sig{0},\sig{1},\dots,\sig{n})[s]=[(\sig{l+1}\dots\sig{n})^{-1}s]$$
Furthermore, $\Lnt{i}$ has a canonical section that is given on
$S$ as $h^*(x_{l+1}\dots x_n)$}
\end{lemma}

\subsection{Local Interpretation of $\Le{i}$}

We can define $i$-admissible charts on $\cm(\cz,\Gamma_Z)$ exactly
as we did on $\cm(\ca,\Gamma_A)$.  Recall that $\Le{i}$ is defined
as $\ev_i^*L$ where $L$ is a line bundle on $Z$ with section $s$
whose zero-scheme is $D$. $\Le{i}$ has a canonical section
$\ev_i^*s$.

\begin{lemma} {Let $S$ be an $i$-admissible chart for
$l\in\{0,1,\dots,n\}$. Then
$$\Lnb{i}|_S\cong 1_S,$$
the topologically trivial bundle under the $\Gt{n}$-action
$$(\sig{1},\dots,\sig{n})[s]=[(\sig{1}\dots\sig{l})^{-1}s]$$
Furthermore, the canonical section is given on $S$ as
$h^*(x_1\dots x_l)$}
\end{lemma}

The proof is exactly analogous to the one given for $\Lnb{i}$.  It
carries through word-for-word with $G[n]$'s substituted for
$\Gt{n}$ and $Z$'s substituted for $A$'s.

\subsection{$\Gamma_1\sqcup \Gamma_2$-admissible charts}

Let $\Gamma$ be a relative graph.  Consider the moduli space
$\cm(\cz,\Gamma)$.

\begin{definition} \label{gadmis} {Let $(\Gamma_Z,\Gamma_A,L,J)$ be
a quadruple whose graph-join, $\Gamma_Z*_{L,J}\Gamma_A$ is
isomorphic to $\Gamma$. Consider a nice chart in
$\cm(\cz,\Gamma)$,
$$\nicechartonZmp{f}{h} $$
Let $$T=S\times_{\cm(\cz,\Gamma)}
\cm(\cz\sqcup\ca,\Gamma_Z\sqcup_{L,J} \Gamma_A).$$ The chart $S$
is said to be {\em $\Gamma_Z\sqcup_{L,J}\Gamma_A$-admissible} if
one of the following happens
\begin{enumerate}
\item[(1)] $T$ is empty.  In this case, $S$ is said to be
trivially admissible.

\item[(2)] There exists an integer $l\in\{1,\dots,n\}$ so that
$T\srarr S$ factors as $T\srarr S\times_{\aff^n} K_l$ so that
$$f:\ex\times_{\aff^n} K_l \srarr Z[n]\times_{\aff^n} K_l = (Z[l-1]\times \aff^{n-l}) \sqcup_D (\aff^{l-1}\times
A[n-l])$$
can be written as
$$g_Z\times h_A\sqcup h_Z\times g_A:\yi_Z\sqcup_{D_l}\yi_A\srarr (Z[l-1]\times \aff^{n-l})\sqcup_D (\aff^{l-1}\times A[n-l])$$
where $\sqcup_{D_l}$ refers to gluing the following families along
$f^{-1}(D_l)$ to form nodes:
$$\xymatrix{
\yi_Z \ar[d] \ar[r]^{f_Z} & Z[l-1]\ar[d]\\
S\times_{\aff^n} K_l \ar@<1ex>[u]^{\gamma_Z,\delta_Z}\ar[r]_{h_Z}
& \aff^{l-1}}$$ and
$$\xymatrix{
\yi_A \ar[d] \ar[r]^{f_A} & Z[n-l]\ar[d]\\
S\times_{\aff^n} K_l
\ar@<1ex>[u]^{\gamma_A,\delta^0,\delta^\infty}\ar[r]_{h_A} &
\aff^{n-l}}$$ where this splitting is described by a quadruple
that is join-equivalent to $$(\Gamma_Z,\Gamma_A,L,J).$$  In this
case, the chart is said to be {\em admissible for $l$}.
\end{enumerate} }
\end{definition}

Again, everything holds for $\cm(\ca,\Gamma)$ if we replace
$(\Gamma_Z,\Gamma_A)$ by $(\Gamma_b,\Gamma_t)$.  By \cite{Li2},
there exists an atlas consisting of
$(\Gamma_Z,\Gamma_A,L,J)$-admissible charts.

\begin{definition} Given a quadruple
$\Upsilon=(\Gamma_Z,\Gamma_A,L,J)$, and a nice chart in
$\cm(\cz,\Gamma)$,
$$\xymatrix{
\ex \ar[d] \ar[r]^{f} & Z[n]\ar[d]\\
S \ar@<1ex>[u]^{\gamma,\delta}\ar[r]_{h} & \aff^n}$$ that is
$\Gamma_Z\sqcup_{L,J}\Gamma_A$-admissible, define $L_\Upsilon$ on
$S$ as follows:
\begin{enumerate}
\item[(1)] If $S$ is trivially admissible, $L_\Upsilon=1_S$, a
topologically trivial bundle with trivial $G[n]$-action.

\item[(2)] If $S$ is admissible for $l\in\{1,\dots,n\}$, define
$L_\Upsilon$ on $S$ as the $G[n]$-equivariant line bundle
$h^*1_\aff^n$ where $1_\aff^n$ is the topologically trivial bundle
with group action
$$(\sig{1},\dots,\sig{n})\cdot s=\sig{l}^{-1}s$$
\end{enumerate}
\end{definition}

Now, $L_\Upsilon$ has a canonical section given by $s_\Upsilon=1$
in trivial charts and $s_\Upsilon=h^*(x_l)$ for non-trivial
charts. It is standard to verify by arguments similar to
\ref{lbiad} that $(L_\Upsilon,s_\Upsilon)$ globalizes to a line
bundle on $\cm(\cz,\Gamma)$ (see \cite{Li2}, Lemma 3.4)

Once we fix $\Gamma$, there are finitely many join-equivalence
classes of quadruples $(\Gamma_Z,\Gamma_A,L,J)$ so that $\Gamma_Z
*_{L,J}\Gamma_A=\Gamma$.  Therefore, we can construct an atlas
that is admissible with respect to every possible join-equivalence
class and is admissible with respect to every interior marked
point.

\subsection{Interpretation of Bundles}

Let us rewrite the bundles $\Dil$,$\Le{i}$,
$\Lnt{i}$,$\Lnb{i}$ as
tensor products of $L_\Upsilon$'s on $\cm(\cz,\Gamma)$ and
$\cm(\ca,\Gamma)$.

On $\cm(\cz,\Gamma)$ where $i$ is the label for an interior marked
point,
\begin{enumerate}
\item[(1)] $\Omega_{\Dil}=\{\Upsilon=(\Gamma_A,\Gamma_Z,L,J)\}$
the set of all join-equivalence classes of quadruples
$\Upsilon=(\Gamma_A,\Gamma_Z,L,J)$.

\item[(2)] $\Omega_{\Le{i}}=\{(\Gamma_A,\Gamma_Z,L,J)|i\in
J(M_A)\}.$
\end{enumerate}

while on $\cm(\ca,\Gamma)$ where $i,j$ are labels for interior
marked points,
\begin{enumerate}
\item[(1)] $\Omega_{\Split}=\{(\Gamma_t,\Gamma_b,L,J)\}.$

\item[(2)] $\Omega_{\Lnb{i}}=\{(\Gamma_t,\Gamma_b,L,J)|i\in
J(M_t)\}$.

\item[(3)] $\Omega_{\Lnt{i}}=\{(\Gamma_t,\Gamma_b,L,J)|i\in
J(M_b)\}$.
\end{enumerate}

\begin{theorem} \label{interpofbundles} {For
$L=\Dil,\Le{i},\Split,\Lnb{i},\Lnt{i}$,
$$L=\bigotimes_{[\Upsilon]\in \Omega_L} L_\Upsilon.$$
where $[\Upsilon]$ denotes a join-equivalence class and $\Upsilon$
a representative element.}
\end{theorem}

\begin{proof} Let us prove (2) on $\cm(\ca,\Gamma)$.  The other cases are similar.  Consider a
chart $S$ that is admissible with respect to every
$\Gamma_b,\Gamma_t$ decomposition and with respect to every
interior marked point:
$$\xymatrix{
\ex \ar[d] \ar[r]^f & Z[n]\ar[d]\\
S \ar@<1ex>[u]^{\gamma}\ar[r]_h & \aff^n}$$
where $\gamma$ is the
section corresponding to the interior marked point $i$.  Suppose
that $S$ is $i$-admissible for $l\in\{0,1,\dots,n\}$. The for any
$k\in\{1,\dots,n\}$ so that $S\times_{\aff^n} K_k$ is not empty,
$$f:\ex\times_{\aff^n} K_k\srarr
A[n]\times_{\aff^n}K_k=(A[k-1]\times\aff^{n-k})\sqcup_D
(\aff^{k-1}\times A[n-k])$$ is given by a join-equivalence class
$[\Upsilon]$ where $\Upsilon=(\Gamma_b,\Gamma_t,L,J)$. Following
from $i$-admissibility, we have the image of every closed point
under
$$(f\circ\gamma):S\srarr (A[k-1]\times\aff^{n-k})\sqcup_D
(\aff^{k-1}\times A[n-k])$$ lying in $\aff^{k-1}\times A[n-k]$ for
$k\leq l$ and lying in $A[k-1]\times\aff^{n-k}$ for $k\geq l+1$.
This implies that $J^{-1}(i)\in M_t$ for $k\leq l$ and
$J^{-1}(i)\in M_b$ for $k\geq l+1$.  Now, let
$B=\{k\in\{1,\dots,l\}|S\times_{\aff^n} K_k=\emptyset\}$ and
define $L_\emptyset$ as a topologically trivial bundle on $\aff^n$
with linearization given by
$$(\sig{0},\sig{1},\dots,\sig{n})s=(\prod_{k\in B}
\sig{k}^{-1})s.$$ $L_\emptyset$ is canonically isomorphic to the
trivial bundle with the trivial linearization.
 Therefore, $$\bigotimes_{[\Upsilon]\in
\Omega_{\Lnb{i}}} L_\Upsilon=\left(\bigotimes_{[\Upsilon]\in
\Omega_{\Lnb{i}}} L_\Upsilon\right)\otimes L_\emptyset=1_S$$ with
the linearization
$$(\sig{0},\sig{1}\dots,\sig{n})\cdot
s=(\sig{1}\dots\sig{l})^{-1}s.$$ But this is just local
interpretation of $\Lnb{i}$ \end{proof}

\subsection{Splitting of Moduli Stacks}

We need to cite a number of results form \cite{Li2}.  These
results were proved for a different moduli stack, $\cm(\cw)$, but
because of the parallels between that space and $\cm(\ca,\Gamma)$
and $\cm(\cz,\Gamma)$, the proofs can be modified in
straightforward fashion.

\begin{definition}\label{multijoin}
Let $\Upsilon=(\Gamma_Z,\Gamma_A,L,J)$ be a quadruples describing
a decomposition in $\cm(\cz,\Gamma*_{L,J}\Gamma_A)$. Define
$m(\Upsilon)$ by
$$m(\Upsilon)=\prod_{i\in RZ} \mu_Z(i).$$
\end{definition}

\begin{theorem} {We have the following equality among cycle classes
$$c_1(L_\Upsilon,s_\Upsilon)\cap\vir{\cm(\cz,\Gamma_Z*_{L,J}\Gamma_A)}=m(\Upsilon)\vir{\cm(\cz\sqcup\ca,\Upsilon)}$$}
\end{theorem}

Let $[\Upsilon]$ be $\Upsilon$'s join-equivalence class.  Then we
have as a consequence of Corollary \ref{totalgluing},

\begin{theorem} {If
$$M_{[\Upsilon]}=\coprod_{(\Gamma_Z',\Gamma_A',L,J)\in[\Upsilon]}
\cm(\ca,\Gamma_A')\times_{D^r}\cm(\cz,\Gamma_Z')$$ is given the
virtual cycle of a disjoint union, then
$$\Psi_{[\Upsilon]}:M_{[\Upsilon]}\srarr \cm(\cz,\Gamma_Z*_{L,J}\Gamma_A)$$
induces
$$\Psi_{[\Upsilon]*}(\vir{M})=|MZ|!|MA|!(|RZ|!)^2\vir{\cm(\cz\sqcup\ca,\Upsilon)}$$}
\end{theorem}

Now, we have the following fiber square for a quadruple
$(\Gamma_Z,\Gamma_A,L,J)$,
$$\xymatrix{\cm(\ca,\Gamma_A)\times_{D^r}
\cm(\cz,\Gamma_Z)\ar[r] \ar[d]
&\cm(\ca,\Gamma_A)\times \cm(\cz,\Gamma_Z)\ar[d]\\
D^r\ar[r]_\Delta & D^r\times D^r}$$ where downward pointing maps
are induced from evaluation at the boundary marked points of
$\cm(\cz,\Gamma_Z)$ and the boundary marked points at $D_\infty$
on $\cm(\ca,\Gamma_A)$ and $\Delta$ is the diagonal.

\begin{theorem} {We have the equality of cycle classes
$$\Delta^!(\vir{\cm(\ca,\Gamma_A)}\times\vir{\cm(\cz,\Gamma_Z)})=\vir{\cm(\ca,\Gamma_A)\times_{D^r}\cm(\cz,\Gamma_Z)}$$}
\end{theorem}

\begin{corollary} If we define
$$P_{\Upsilon}=\coprod_{(\Gamma_Z',\Gamma_A',L,J)\in[\Upsilon]}
\cm(\ca,\Gamma_A')\times \cm(\cz,\Gamma_Z')$$ then
$$c_1(L_\Upsilon,s_\Upsilon)\cap\vir{\cm(\cz,\Gamma_Z*_{L,J}\Gamma_A)}=\frac{m(\Upsilon)}{|MZ|!|MA|!(|RZ!)^2}\Delta^!(P_\Upsilon)$$
\end{corollary}
$L$ together with $i:X\srarr Z$ induces a morphism
$$\Lambda:(X^{|MA|}\times X^{|RA_0|}) \times Z^{|MZ|}\srarr
Z^M\times X^R$$ where $M=|MZ|+|MA|$ and $R=|RA_0|$ are the number
of interior and boundary marked points in $\Gamma_A *_{L,J}
\Gamma_Z$

We have morphisms
$$\xymatrix{
X^{|MA|+|RA_0|}\times X^{|RZ|} \times Z^{|MZ|}
\ar[r]^<<<<<{\tilde{\Delta}} \ar[d]^p & X^{|MA|+|RA_0|} \times
X^{|RA_\infty|}
\times Z^{|MZ|} \times X^{|RZ|}\\
(X^{|MA|}\times X^{|RA_0|}) \times Z^{|MZ|}&& }$$ where
$\tilde{\Delta}$ is induced by $\Delta:X^{|RZ|}\srarr
X^{|RA_\infty|}\times X^{|RZ|}$ and $p$ is the projection.

Therefore, for $c\in H^*(Z^{|MZ|}\times X^{|RZ|}$,
$$\deg((\Ev^*(c)\cup
c_1(L_\Upsilon))\cap\vir{\cm(\cz,\Gamma_A*_{L,J}\Gamma_Z)})$$
$$=\frac{m(\Upsilon)}{\Aut_{\Gamma_Z,\Gamma_A,L}(RZ,RA_\infty)}
\deg(\Ev^*({\tilde{\Delta}}_!(p^*\Lambda^*c)\cap(\vir{\cm(\ca,\Gamma_A)}\times\vir{\cm(\cz,\Gamma_Z)})$$




Again, there are obvious degeneration formulae on
$\cm(\ca,\Gamma_b*_{L,J}\Gamma_t)$ obtained by replacing
$\Gamma_Z$ with $\Gamma_b$ and $\Gamma_A$ with $\Gamma_t$.

By writing
$$c_1(L_\Omega)=\sum_{[\Upsilon]\in\Omega} c_1(L_\Upsilon)$$
we obtain formulas for expressing
$$\deg((\Ev^*(C)\cup
c_1(L_\Omega))\cap\vir{\cm(\cz,\Gamma_Z*_{L,J}\Gamma_A)})$$ and
$$\deg((\Ev^*(C)\cup
c_1(L_\Omega))\cap\vir{\cm(\ca,\Gamma_b*_{L,J}\Gamma_t)}).$$

\subsection{Normal Bundle to Split Maps}

The split maps in $\cm(\cz,\Gamma)$ for a divisor in
$\cm(\cz,\Gamma)$. It is natural to ask what the normal bundle to
such a divisor is. Since a localization computation will have
split maps as fixed loci, knowledge of the normal bundle will be
indispensable for localization as in \cite{GV}. This section is
logically independent from the rest of this paper, but is included
as reference.

Consider $(Z,D)$, a projective manifold $Z$ and a smooth divisor
$D\subset Z$.  Let $(X,L)$ be given by $X=D$, $L=N_{D/Z}$.  Let
$\Gamma_Z$ and $\Gamma_A$ be relative and rubber graphs,
respectively. Let $L,J$ be joining data as in \ref{graphjoin}.
Consider $L_\Upsilon$ for the quadruple
$\Upsilon=(\Gamma_Z,\Gamma_A,L,J)$. Then $c_1(L_\Upsilon)$ is a
substack of $\cm(\cz,\Gamma_Z*_{L,J}\Gamma_A)$.

Consider the moduli stacks $\cm(\cz,\Gamma_Z)$,
$\cm(\ca,\Gamma_A)$, $\cm(\cz,\Gamma_Z*_{L,J}\Gamma_A)$, and the
inclusion
$$\Phi:\cm(\ca,\Gamma_A)\times_{D^r} \cm(\cz,\Gamma_Z)\srarr \cm(\cz,\Gamma_Z*_{L,J}\Gamma_A).$$
$\cm(\ca,\Gamma_A)\times_{D^r} \cm(\cz,\Gamma_Z)$ has projections
$p_A,p_Z$ to its factors.

\begin{proposition} {$\Phi^*L_\Upsilon=p_Z^*\Dil^\vee \otimes
p_A^*\Bot^\vee$.}
\end{proposition}

\begin{proof} Consider an $\Upsilon$ admissible chart for
$l\in\{1,\dots,n\}$,
$$\xymatrix{
\ex \ar[d] \ar[r]^f & Z[n]\ar[d]\\
S \ar@<1ex>[u]^{\gamma,\delta}\ar[r]_h & \aff^n}$$
Then
$L_\Upsilon$ is given as $h^*1_{\aff^n}$ where $1_{\aff^n}$ is the
topologically trivial bundle with group action
$$(\sig{1},\dots,\sig{n})\cdot s=(\sig{l})^{-1}s.$$
Let $K_l\subseteq\aff^n$ be the subset given by $t_l=0$.  Let
$\phi:K_l\srarr \aff^n$ be the inclusion.

\begin{lemma} There is an isomorphism of swags, $K_l\cong
\aff^{l-1}\times \aff^{n-l}$ where $\aff^{l-1}$ is considered a
swag under $G[l-1]$, $\aff^{n-l}$ is considered a swag under
$\Gt{n-l}$.  Under this isomorphism,
$$\phi^*L_\Upsilon\cong \pi_1^* DIL[l-1]^\vee \otimes \pi_2^* B[n-l]^\vee$$
where $\pi_1$ and $\pi_2$ are defined as
$$\xymatrix{
&K_l\ar[dl]_{\pi_1}\ar[dr]^{\pi_2}&\\
\aff^{l-1} & & \aff^{n-l}}$$
\end{lemma}

\begin{proof} The isomorphism is given by
\begin{eqnarray*}
H_l&\srarr&\aff^{l-1}\times\aff^{n-l}\\
(t_1,\dots,t_{l-1},0,t_{l+1},\dots,t_n)&\mapsto&
(t_1,\dots,t_{l-1}),(t_{l+1},\dots,t_n)\\
G[n]&\srarr &G[l-1]\times \Gt{n-l}\\
(\sig{1},\dots,\sig{n})&\mapsto&
(\sig{1},\dots,\sig{l-1}),(\sig{1}\sig{2}\dots\sig{l},\sig{l+1},\dots,\sig{n})
\end{eqnarray*}

Note that in the above, $\sig{1}\dots\sig{l}$ gets mapped into the
zeroth place in $\Gt{n-l}$.

Therefore, $\pi_1^* DIL[l-1]^\vee\otimes \pi_2^* B[n-l]^\vee$, is
the topologically trivial bundle on $H_l$ with the group action
$$(\sig{1},\dots,\sig{n})\cdot
s=(\sig{1}\dots\sig{l-1})(\sig{1}\dots\sig{l})^{-1}s=\sig{l}^{-1}s$$
\end{proof}

The above isomorphism is canonical and globalizes giving the
conclusion. \end{proof}



    \section{The Trivial Cylinder Theorem}

\subsection{Trivial Cylinders}

We will single out certain connected components of curves
parameterized by $\cm(\ca,\Gamma)$.  These are the {\em trivial
cylinders} which will be significant in \cite{Ka2}.

\begin{definition} \label{trivialcyl} Let $\Gamma$ be a rubber
graph. A vertex $v$ is said to correspond to a {\em trivial
cylinder} if
\begin{enumerate}
\item[(1)] $g(v)=0$.

\item[(2)] $d(v)=0$

\item[(3)] $a_0^{-1}(v)$ is a single point.

\item[(4)] $a_\infty^{-1}(v)$ is a single point.

\item[(5)] $\mu^0(a_0^{-1}(v))=\mu^\infty(a_\infty^{-1}(v))$.

\item[(6)] $A_M^{-1}(v)$ is empty.
\end{enumerate}
\end{definition}

A trivial cylinder corresponds to a component $\ex_v$ in a family.
This component is given over a point $p\in S$ by a map $f$ from
$C$ to $\ _l P$ where $C$ is a chain of $l$ $\proj^1$'s.  $f$
takes $C$ into a fiber over some point $x\in X$, and there is an
integer $d$ so that $f$ is given on each $\proj^1$ by
$$z\mapsto z^d.$$
Note that if $\Gamma$ consists of a single vertex corresponding to
a trivial cylinder, then a morphism of a family of type $\Gamma$,

$$\xymatrix{
\ex \ar[d] \ar[r]^f & A[n]\ar[d]\\
S \ar@<1ex>[u]^{\delta^0_1,\delta^\infty_1}\ar[r] & \aff^n}$$
would be invariant under composing $f$ with the action of
$(\sig{0},1,\dots,1)\in\Gt{n}$ and pre-composing with the
appropriate $S$-isomorphism on $\ex$. Consequently a trivial
cylinder is not stable. This does not rule out morphisms of type
$\Gamma$ which has a component which is a trivial cylinder.

There is a straightforward way of comparing relative Gromov-Witten
invariants involving graphs with trivial cylinders to those
without.  It follows from the following theorem.

\begin{theorem}\label{tctheorem} {Let $\Gamma$ be some rubber graph.  Let $\Gamma_|$ be $\Gamma$
together with a vertex corresponding to a trivial cylinder.  There
is a natural map
$$v:\cm(\ca,\Gamma_|)\srarr \cm(\ca,\Gamma)\times X$$
so that
$$v_*\vir{\cm(\ca,\Gamma_|)}=\frac{1}{r}\vir{\cm(\ca,\Gamma)}\times
[X]$$ and
$$v^*(\Bot)=\Bot.$$
Consequently if $\Gamma$ has $m$ interior marked points and
$r_0+r_\infty$ boundary marked points then we have a commutative
diagram
$$\xymatrix{
\cm(\ca,\Gamma_|)\ar[d]\ar[r]^>>>>>>>>>{\Ev_|}& X^m\times
X^{r_0+1}\times
X^{r_\infty+1}\ar[d]^h\\
\cm(\ca,\Gamma)\times X\ar[r]_>>>>>{\Ev\times\Delta}&X^m\times
X^{r_0}\times X^{r_\infty}\times (X\times X)}$$ where
$\Delta:X\srarr X^2$ is the diagonal and the morphism $h$ reorders
the products of $X$ so that the copies of $X$ corresponding to the
$r_0+1$st and $r_\infty+1$st boundary marked points are taken to
the image of the diagonal morphism. Then for classes $A\in
H^*(X^m\times X^{r_0+1}\times X^{r_\infty+1})$, we have
$$\deg(\Ev_|^*(A)\cap\vir{\cm(\ca,\Gamma_|)})=\frac{1}{r}\deg(((\Ev\times\Delta)\circ h^{-1})^*(A)\cap
(\vir{\cm(\ca,\Gamma)}\times [X]))$$ }
\end{theorem}

There is a heuristic argument for why this is true.  We can think
of a trivial cylinder in $A[0]$ as a map $\proj^1\srarr P$ given
by
$$z\mapsto z^r$$
Note that this map $r$ automorphisms given by
$$z\mapsto \omega_r^a z$$
where $\omega_r$ is an $r$th root of unity and $0\leq a\leq r-1$.
Therefore, we would like to say
$$\cm(\ca,\Gamma_|)=\cm(\ca,\Gamma)\times (X/(\zee/r))$$
where $X/(\zee/r)$ is the stack consisting of $X$ quotiented by a
cyclic group of order $r$, acting trivially.  We should even hope
for the virtual cycles of each space to agree.  This,
unfortunately, is not true as stated.

The main obstacle for this fact being the case is that the target
of a map in $\cm(\ca,\Gamma)$ may not be $P$, but $_l P$ for
$l\geq 1$.  Therefore, we may have split maps   This gives an
automorphism group of $(\zee/r)^{l+1}$ where we can multiply by a
different power of $\omega_r$ on each $P$.  Moreover, this
splitting phenomenon gives a non-reduced scheme where points
corresponding to maps into a chain of length $l$ have multiplicity
$r^l$. Therefore, we'd like to define a map
$$v:\cm(\ca,\Gamma_|)\srarr \cm(\ca,\Gamma)\times X$$

If we have
$$v_*(\vir{\cm(\ca,\Gamma_|)})=\frac{r^l}{r^{l+1}}\vir{\cm(\ca,\Gamma)}\times X$$

then we can conclude the theorem.

The argument below was suggested by Jun Li.

\subsection{Stacks of Trivial Cylinders}

We need to define a stack $\mt$ that parameterizes map rubber maps
of cylinders.  Since trivial cylinders are not stable, this moduli
stack is not Deligne-Mumford.

Let $\Delta_r$ be a rubber graph corresponding to a degree $r$
trivial cylinder.  Consider $\mt=\cm(\ca,\Delta_r)$ constructed as
above but where we do not impose the stability condition on
families. Instead, we just impose that for a family over any
point,
$$\xymatrix{C\ar[r]^f\ar[d]&A[n]\ar[d]\\
p\ar[r]&\aff^n}$$
if $A[n]\times_{\aff^n} p$ is a chain of $P$'s,
$$P_0\sqcup_D P_1\sqcup_D\dots \sqcup_D P_l$$ then $C$ consists of a chain
$$C_0\cup C_1\cup\dots\cup C_l$$ where $f$ maps $C_i$ to a
fiber of $\pi:P_i\srarr X$ by a map of the form
$$z\srarr z^r$$
There is a natural (although not representable) map
$v:\mt\srarr\ca^{\rel}\times X$ taking each family of trivial
cylinders to its family of target schemes together with the point
in $X$ in whose fiber it lies.

We can get an explicit handle on the morphism $v$ by pulling it
back by the morphism $\aff^n\srarr \ca^{\rel}$ given by the family
of targets $A[n]\srarr\aff^n$.  Now $T=\mt\times_{\ca^{\rel}}
\aff^n$ has an explicit description. Define the scheme $\TC^n$ by
$$\TC^n=\Spec
\kay[x_1,\dots,x_n,y_1,\dots,y_n]/(x-y_1^r,\dots,x-y_n^r).$$
$TC^n\times X$ has a natural morphism to $T$.  The map to $\aff^n$
is given by
$$TC^n\srarr \aff^n=\Spec \kay[x_1,\dots,x_n].$$
The map to $\mt$ comes from the following family: the domain is
given by
$$A[n]\times_{\aff^n\times X} TC^n$$
where the map $TC^n\srarr\aff^n$ is
$$TC^n\srarr\aff^n=\Spec \kay[y_1,\dots,y_n]$$
and $A[n]\srarr\aff^n\times X$ is the usual projection $p\times
\pi$. The target is given by
$$A[n]\times_{\aff^n\times X} TC^n$$
where the map $TC^n\srarr\aff^n$ is
$$TC^n\srarr\aff^n=\Spec \kay[x_1,\dots,x_n].$$
The stable morphism of families, $f$ is given by constructing
nodal normal neighorhoods $w_1w_2=x_1$ and $z_1z_2=y_1$ and
setting $$f^*w_i=z_i^r.$$ This induces a degree $r$ covering of
the fibers.

This morphism has the automorphisms
$$TC^n\times_T \TC^n=(\zee/r)^{n+1}\times TC^n.$$
The first projection is
$$(a_0,a_1,\dots,a_n)\cdot (x_1,\dots,x_n,y_1,\dots,y_n)\mapsto
(x_1,\dots,x_n,y_1,\dots,y_n)$$ where $a_i\in\zee/r$.  The second
projection is
$$(a_0,a_1,\dots,a_n)\cdot (x_1,\dots,x_n,y_1,\dots,y_n)\mapsto
(x_1,\dots,x_n,\omega_r^{a_1}y_1,\dots,\omega_r^{a_n}y_n)$$

The fact that $y$ becomes nilpotent at $x=0$ come from the fact
that when we consider a split curve, we must smooth the node to
$r$th order before we can move the curve out of its singular
target fiber.  The $\zee/r$ factors come from rotating different
components of split maps.

Consider the map
$$v:\mt\times_{\ca^\rel} \aff^n\srarr
(\aff^n\times_{\ca^\rel}\ca^\rel)\times X=\aff^n\times X.$$ It is
clear that for any point $q\in\aff^n\times X$, we have
$\deg(v^{-1}(q))=\frac{1}{r}$. This follows from the fact that
degree and local multiplicity are properties that can be checked
at generic geometric points.

Note that given any map $\Spec \kay\srarr \ca^\rel$, we may find
an $n$ so that the map factors as $\Spec
\kay\srarr\aff^n\srarr\ca^\rel$.  Consequently, $\deg(\Spec
\kay\times_{\cz^\rel}\mt)=\frac{1}{r}$.

\subsection{Comparing $\cm(\ca,\Gamma_|)$ to $\cm(\ca,\Gamma)$}

\begin{theorem}
{$\cm(\ca,\Gamma_|)=\cm(\ca,\Gamma)\times_{\ca^\rel} \mt$}
\end{theorem}

\begin{proof} This is a matter of unwinding the definition of a fiber
product in the category of stacks.  Let
$p_\ma,p_T:\cm(\ca,\Gamma),\mt\srarr \ca^{rel}$ be the
projections. A family $S\srarr\cm(\ca,\Gamma)\times_{\ca^\rel}\mt$
consists of a triple $(h_\ma,h_T,a)$ where
$$h_\ma:S\srarr \cm(\ca,\Gamma)$$
$$h_T:S\srarr \mt$$
and $a$ is an arrow in $\ca^{\rel}$, $a:p_\ca(h_\ma)\srarr
p_T(h_T)$.  If we shrink $S$ to make $h_\ca$ and $h_T$ a nice
family, we get diagrams
$$\xymatrix{\ex\ar[r]^<<<<<<{f_\bullet}\ar[d]&A[n_1]\ar[d]\\
S\ar[r]&\aff^{n_1}}$$
and
$$\xymatrix{T\ar[r]^<<<<<<{f_T}\ar[d]&A[n_2]\ar[d]\\
S\ar[r]&\aff^{n_2}}$$
By shrinking $S$ further, we can find an
effective arrow between $S\srarr \aff^{n_1}$, $S\srarr
\aff^{n_2}$, that is an integer $N$, standard inclusions
$[n_1]\hookrightarrow [N]$, $[n_2]\hookrightarrow [N]$ and maps
$\rho_1,\rho_2:S\srarr \Gt{N}$ so that the following compositions
are equal
$$S\srarr S\times\aff^{n_1}\srarr S\times\aff^N \srarr
\Gt{N}\times\aff^N\srarr\aff^N$$
$$S\srarr S\times\aff^{n_2}\srarr S\times\aff^N \srarr
\Gt{N}\times\aff^N\srarr\aff^N.$$ This allows us to combine the
two families
$$\xymatrix{
\ex\sqcup T\ar[r]^>>>>>{f_| \sqcup f_T}\ar[d] & A[N]\ar[d] \\
S\ar[r]&\aff^{N}}$$ to obtain a family in $\cm(\ca,\Gamma_|)$.
Showing that this map from
$\cm(\ca,\Gamma)\times_{\ca^\rel}\mt\srarr \cm(\ca,\Gamma_|)$ is
an isomorphism is a straightforward verification. \end{proof}

\begin{corollary} {There is a natural morphism
$\cm(\ca,\Gamma_|)\srarr \cm(\ca,\Gamma)\times X$}
\end{corollary}

\begin{proof} Consider $v:\mt\srarr\ca^\rel\times X$.  This induces
$$v:\cm(\ca,\Gamma_|)=\cm(\ca,\Gamma)\times_{\ca^\rel} \mt\srarr
\cm(\ca,\Gamma)\times_{\ca^\rel} (\ca^\rel\times
X)=\cm(\ca,\Gamma)\times X.$$ \end{proof}

\subsection{Virtual Cycles}

In this section we will show how to compare the virtual cycle of
$\cm(\ca,\Gamma_|)$ to that of $\cm(\ca,\Gamma)\times X$.

Let us first note that $\mt$ has no obstructions.  Let
$V\srarr\ca^\rel$ be some family.  Then consider the
Deligne-Mumford stack $\mt\times_{\ca^\rel} V$.  We can write down
the tangent/obstruction complex for $\mt\times_{\ca^\rel} V$
following the recipe in \cite{Li2}.  We see that
$h^2(E^\bullet)=0$, so
$$\vir{\mt\times_{\ca^\rel} V}=\mt\times_{\ca^\rel} V$$

Now let us look at the stack $\cm(\ca,\Gamma)\times_{\ca^\rel}
\mt$.  We consider a particular kind of \'{e}tale chart.  Let
$S\srarr\cm(\ca,\Gamma)$ be \'{e}tale.  Then let $W\srarr
S\times_{\ca^\rel} \mt$ be an \'{e}tale chart.  Then we have the
following composition of \'{e}tale morphisms
$$W=S\times_S W\srarr S\times_S
(S\times_{\ca^\rel}\mt)=S\times_{\ca^\rel}\mt\srarr
\cm(\ca,\Gamma)\times_{\ca^\rel} \mt$$

If we shrink $W$ as above, we can obtain nice families
$$\xymatrix{\ex\ar[r]^{f_\bullet}\ar[d]&A[N]\ar[d]\\
W\ar[r]_h&\aff^N}$$ and
$$\xymatrix{T\ar[r]^{f_T}\ar[d]&A[N]\ar[d]\\
S\ar[r]&\aff^N}$$
Now, the tangent/obstruction complex on
$W=S\times_S W$ splits into contributions from
$S\srarr\cm(\ca,\Gamma)$ and $W\srarr\mt$. Moreover, the induced
tangent/obstruction complex is identical to the one given by
considering $W\srarr\cm(\ca,\Gamma_|)$. Therefore, we obtain a
cone cycle $C\in Z_*(E=E_\ca\oplus E_T)$ over $W$ where
$$C=C_\ca\times_S C_T$$
where $C_\ca\in Z_*(E_\ca)$ is the cone cycle on $S$ and $C_T\in
Z_*(E_T)$ on $W$.  By the naturality of the Gysin map (see
\cite{Fulton} 6.5),
$$s_{E_\ca\oplus E_T}^!(C_\ca\times_S
C_T)=(s_{E_\ca}^! C_\ca)\times_S W$$ Globally, this construction
tells us that
$$\vir{\cm(\ca,\Gamma_|)}=\vir{\cm(\ca,\Gamma)\times_{\ca^\rel}\mt}=
\vir{\cm(\ca,\Gamma)}\times_{\ca^\rel}[\mt]$$

Therefore, if we consider the morphism,
$$v:\cm(\ca,\Gamma_|)\srarr \cm(\ca,\Gamma)\times X$$ induced
from the degree $\frac{1}{r}$ morphism
$$\mt\srarr\ca^\rel\times X,$$
we have
$$v_*(\vir{\cm(\ca,\Gamma_|)})=v_*(\vir{\cm(\ca,\Gamma)}\times_{\ca^\rel}[\mt])=
\frac{1}{r}\vir{\cm(\ca,\Gamma)}\times [X].$$


    \bibliographystyle{plain}
    \bibliography{lbosorm}

\end{document}